\documentclass[preprint]{elsarticle}
\usepackage[ansinew]{inputenc}
\usepackage{a4wide}
\usepackage{calc}
\usepackage{times}
\usepackage{amssymb,amsmath,amstext,amsthm,amsfonts, ams fonts}
\usepackage{colortbl,color,graphics,graphicx,epsfig}
\definecolor{gray}{rgb}{0.8,0.8,0.8}
\usepackage{dsfont}
\usepackage{bbm}
\usepackage{nicefrac}
\usepackage[pdflinkmargin=5pt,pdfstartview={FitBH -32768},plainpages=false]{hyperref}
\newcommand{\dunion}{\dot\bigcup}

\renewcommand{\P}{\mathbb{P}}
\newcommand{\F}{\mathcal{F}}
\renewcommand{\H}{\mathcal{H}}
\newcommand{\G}{\mathcal{G}}
\newcommand{\E}{\mathbb{E}}
\newcommand{\N}{\mathds{N}}
\newcommand{\R}{\mathds{R}}
\newcommand{\1}{\mathbbm{1}}
\renewcommand{\d}{\Delta}
\newcommand{\KLEINO}{{\scriptstyle{\mathcal{O}}}}

\DeclareMathAccent{\verywidehat}{\mathord}{largesymbols}{'144}
\newcommand{\var}{\mathbb{V}\hspace*{-0.05cm}\textnormal{a\hspace*{0.02cm}r}}

\DeclareMathOperator{\AVAR}{\mathbf{AVAR}}
\renewcommand{\:}{\mathrel{\mathop{:}}}

\newdefinition{remark}{Remark}
\newdefinition{defi}{Definition}
\newtheorem{theo}{Theorem}
\newtheorem{assump}{Assumption}
\newtheorem{prop}{Proposition}[section]
\newtheorem{lem}[prop]{Lemma}
\newtheorem{cor}[prop]{Corollary}

\begin{document}
\begin{frontmatter}



\title{Asymptotics of Asynchronicity}

\author{Markus Bibinger\fnref{label1}}
\fntext[label1]{Financial support from the Deutsche Forschungsgemeinschaft via SFB 649 `\"Okonomisches Risiko', Humboldt-Universität zu Berlin, is gratefully acknowledged.}
\address{Institut f\"ur Mathematik, Humboldt-Universit\"at zu Berlin, Unter den Linden 6, 10099 Berlin, Germany}

\begin{abstract}
In this article we focus on estimating the quadratic covariation of continuous semimartingales from discrete observations that take place at asynchronous observation times. The Hayashi-Yoshida estimator serves as synchronized realized covolatility for that we give our own distinct illustration based on an iterative synchronization algorithm. We consider high-frequency asymptotics and prove a feasible stable central limit theorem. The characteristics of non-synchronous observation schemes affecting the asymptotic variance are captured by a notion of asymptotic covariations of times. These are precisely illuminated and explicitly deduced for the important case of independent time-homogeneous Poisson sampling.
\end{abstract}

\begin{keyword}
non-synchronous observations \sep quadratic covariation \sep Hayashi-Yoshida estimator \sep stable limit theorem\sep asymptotic distribution\\ \vspace*{.3cm}\noindent
\textit{MSC Classification:} 62M10\sep 62G05\sep 62G20\sep 91B84\\ \vspace*{.3cm}\noindent
\textit{JEL Classification:} C14\sep C32\sep C58\sep G10
\end{keyword}
\end{frontmatter}
\thispagestyle{plain}
\section{Introduction\label{sec:1}}
Nonparametric estimation methods for the quadratic variation of semimartingales have become an issue of great interest in recent years. One reason is the interpretation of the quadratic variation of the continuous part as integrated volatility in financial modeling.\\ If a semimartingale is observed discretely at times $t_i,0\le i\le n$ on a finite time horizon $[0,T]$, the sum of squared returns (increments of the semimartingale), called realized volatility, converges to the quadratic variation as $\sup{\left(t_i-t_{i-1}\right)}\rightarrow 0$ as $n\rightarrow\infty$. The same fact pertains to the multi-dimensional case where the realized covolatilities of two processes converge to the quadratic covariations. More usually multivariate data, in particular financial time series, are recorded at times following non-synchronous observation schemes. Therefore, realized covolatility estimates most commonly incorporate a previous-tick interpolation approach. Though, this machinery leads to the so-called Epps effect \cite{epps} that realized covolatilities tend to zero as the sampling frequency increases. Especially for the more and more available ultra high-frequency financial tick-data this issue poses problems.\\
A solution for the asynchronous estimation problem has been proposed in \cite{hy}. We call this estimator which arises as realized covolatility from all products of returns with overlapping observation time instants Hayashi-Yoshida estimator. Our investigation of that estimation approach leads to several useful rewritings and interpretations. The final representation is based on an iterative synchronization procedure which has been used first in \cite{palandri}. This synchronized realized covolatility and the data aggregation technique for synchronization can serve as a basis for combined approaches in various generalizations of the underlying statistical model. A very important enhancement of the model in that we take market microstructure noise into account is covered in \cite{bibinger} and \cite{bibinger2} by extending the synchronized realized covolatility to a generalized multiscale estimator. \\
The asymptotic theory developed in the article on hand is grounded on stable limit theorems for semimartingales from \cite{jacod1}. We obtain a stable limit theorem for the process associated with the estimation error of a Hayashi-Yoshida estimator for the quadratic covariation at time $t\in[0,T]$ and for the overall estimator by the marginal distribution at $t=T$ the stable weak convergence to a centred mixed normal limiting distribution. The random asymptotic variance splits up in two terms induced by an idealized synchronous approximation and an additional error due to the lack of synchronicity.\\
The article is arranged in six sections. In the following Section \ref{sec:2} we give insight into the concept of stable weak convergence and a short review on the essential theory from the literature. In Section \ref{sec:3} the Hayashi-Yoshida estimator and our related synchronization algorithm that we have first presented in \cite{bibinger} is revisited and the asymptotic theory including the key result is provided in Section \ref{sec:4}. The detailed proof of the central Theorem \ref{HYclt} is postponed to the Appendix \ref{sec:7}. Following some simple illustrative and motivating examples before, Section \ref{sec:5} comes up with the analysis for the important time-homogeneous independent Poisson sampling case for which we evaluate all ingredients of the asymptotic variance explicitly. To benefit from the stable central limit theorem and provide a basis for statistical inference, we give a consistent estimator for the asymptotic variance in Section \ref{sec:6}. 
\section{Stable convergence and Jacod's stable limit theorem revisited\label{sec:2}}
This section is devoted to the notion of stable weak convergence which will be an essential concept for the development of our limit theory throughout this article. The concept of stable convergence goes back to \cite{renyi} and results about stable limit theorems were extended in \cite{aldous} and \cite{feigin}. The reason what makes stable weak convergence a key element of our asymptotic considerations, is that it allows to conclude joint weak convergence when we derive results about asymptotic mixed normality. In the case that a sequence of random variables $(X_n)$ weakly converges to a mixed Gaussian limiting random variable $VZ$, with $Z$ being standard normally distributed, $Z\sim\mathbf{N(0,1)}$, and a strictly positive random variable $V$, independent of $Z$, we cannot derive confidence intervals if the distribution of $V$ is unknown. However, if a consistent estimator $V_n^2$ for the asymptotic variance $V^2$ is available (in the sense that $V_n^2\stackrel{p}{\longrightarrow}V^2$), the stable weak convergence will assure that $(X_n,V_n^2)\rightsquigarrow(VZ,V^2)$ jointly and also that $X_n/V_n\rightsquigarrow Z$. The last implication also holds in a stable version: $X_n/V_n\stackrel{st}{\rightsquigarrow} Z$. For that reason, if we are in situations as described above we gain from proving stable weak convergence which paves the way towards statistical inference.\\
Next, we present the formal definition and the main properties of stable weak convergence of sequences of random variables.\\ \\
Let $(X_n)$ be a sequence of random variables defined on some probability space $(\Omega,\mathcal{F},\P)$ and taking values in a Polish space $(E,\mathcal{E})$. We say that the sequence $(X_n)$ converges weakly in $\mathbb{L}^1$ to $X$ if for any bounded random variable $Z$
$$\lim_{n\rightarrow\infty}\E\left[ZX_n\right]=E\left[ZX\right]$$
holds.\\
\begin{defi}\label{stable}For a sub-$\sigma$-field $\mathcal{G}\subseteq\mathcal{F}$ the sequence of random variables $(X_n)$ is said to converge $\mathcal{G}$-stably, if there is a random probability measure $\mu$ on $(\Omega\times E,\mathcal{G}\otimes\mathcal{E})$ such that
$$\lim_{n\rightarrow\infty}\E\left[Zf(X_n)\right]=\int_{\Omega\times E}\mu(d\omega,dx)Z(\omega)f(x)$$
for all $f\in\mathcal{C}_b(E)$ (continuous and bounded) and $\mathcal{G}$-measurable bounded random variables $Z$.\\
If $\mathcal{G}=\mathcal{F}$, we say $(X_n)$ converges stably in law to $X$ ($X_n\stackrel{st}{\rightsquigarrow}X$).
\end{defi}
\begin{remark}$\mathcal{G}$-stable convergence is the weak convergence in $\mathbb{L}^1$ of $\E\left[f(X_n)|\mathcal{G}\right]$ for all $f\in\mathcal{C}(E)$ to $\mu\circ f$. This implies convergence in distribution to the probability measure $\nu$ defined by 
$$\nu(B)=\int\mu(d\omega,B)\1_{\{X(\omega)\in B\}}\P(d\omega)~.$$
If $(X_n)$ converges stably, the limiting law is $\mu(\Omega,\,\cdot\,)$.
\end{remark}\noindent
The following proposition states some useful equivalent characterizations of $(\mathcal{G}-)$stable convergence.
\begin{prop}$(X_n)$ converges $\mathcal{G}$-stably is equivalent to:
\begin{enumerate}
\renewcommand{\labelenumi}{(\roman{enumi})}
\item For every $\mathcal{G}$-measurable random variable $Z$ on $\Omega$, $(Z,X_n)$ converges in law.
\item For every $\mathcal{G}$-measurable random variable $Z$ on $\Omega$, $(Z,X_n)$ converges $\mathcal{G}$-stably.
\item The sequence $(X_n)$ is tight, and for all $G\in\mathcal{G}$ and $f\in\mathcal{C}(E)$, the sequence $\E\left[\1_{G}f(X_n)\right]$ converges.
\end{enumerate}
\end{prop}
This proposition is proved in \cite{shir} as part of Proposition IX.1.4.
Stable (weak) convergence is a stronger mode of ordinary convergence in distribution. It is weaker than convergence in probability, but we emphasize that the limit depends on the limiting random variable $X$ itself and not only on the distribution of $X$.\\
If $(X_n)$ converges stably to $X$, $X$ is defined on an extension $(\Omega^{\prime},\mathcal{F}^{\prime},\P^{\prime}=\mu)$ of the original probability space, so that
$$\forall f\in\mathcal{C}(E):~~\lim_{n\rightarrow\infty}\E\left[f(X_n)Z\right]=\E^{\prime}\left[f(X)Z\right]~.$$
In the situation that we face in this article, a Gaussian random variable which is independent of $\mathcal{F}$ will appear as limiting law. In this case we call the extension of the original probability space orthogonal. In particular, if $X_n\stackrel{st}{\rightsquigarrow}X$ holds with an $\mathcal{F}$-measurable random variable $X$ ($\mathbb{L}^1-$convergence), the foregoing proposition yields that $(X_n,X)\rightsquigarrow(X,X)$ and, hence $(X_n-X)\rightsquigarrow 0$ holds, which implies convergence in probability. Therefore, in all cases where stable weak convergence is a suitable adequate concept, limiting laws are defined on a genuine extension of the original probability space.\\
The following proposition gives the result that stable convergence is the suitable concept to derive feasible central limit theorems and confidence intervals if the asymptotic variances in limit theorems are unknown random, but can be estimated consistently.
\begin{prop}Let $(X_n,V_n)$ be real-valued random variables defined on $(\Omega,\mathcal{F},\P)$. If $X_n\stackrel{st}{\rightsquigarrow} X$ with a mixed normal limiting random variable $X\sim \mathbf{N}(0,V^2)$ and $V_n\stackrel{p}{\longrightarrow}V$ with $V$ being $\mathcal{F}$-measurable. Then
$$X_n/V_n\stackrel{st}{\rightsquigarrow}\mathbf{N}(0,1)$$
holds true.\end{prop}
Note that we use the same denotation expression for mixed normal laws and common normal laws and the difference becomes clear out of the context and by the specific variances.
On the assumptions of the proposition $(X_n,V_n)\stackrel{st}{\rightsquigarrow}(X,V)$ is implied and the convergence of $X_n/V_n$ follows by the continuous mapping theorem. This proposition is part of Proposition 2.5 in \cite{poldilimit}. We restricted ourselves to real-valued random variables in the last proposition. A more general version can be found in \cite{shir}. \\
The concept of stable convergence also carries over to stochastic processes. For this extension of stable convergence to stochastic processes, or more precisely to semimartingales, the Polish space $E$ in Definition \ref{stable} is chosen to be the Skorohod space.
The following limit theorem for stable convergence of continuous local martingales will be the foundation for our later deduced limit theorem in this article: 
\begin{theo}[\textbf{Jacod's theorem: A martingale version}]
\label{jacod}If $(M_t,\mathcal{F}_t)$ with $0\le t<\infty$ is a continuous local martingale defined on the probability space $(\Omega,\mathcal{F},\P)$, we denote by $\mathcal{M}^{\perp}$ the set of bounded $(\mathcal{F}_t)$-adapted martingales orthogonal to $M=(M_t,\mathcal{F}_t)$ what means that $\left[ M,M^{\perp}\right]\equiv 0$. If $(X^n)$ is a sequence of continuous $(\mathcal{F}_t)$-adapted local martingales for which
\begin{equation}\left[ X^n\right]_t\stackrel{p}{\longrightarrow}V_t~~\forall t\,\in\,[0,\infty)\end{equation}
with a continuous process $V$ holds, the following two conditions 
\begin{subequations}
\begin{equation}\label{conditionsstable1}\left[ X^n,M\right]_t\stackrel{p}{\longrightarrow}0~~\forall t\,\in\,[0,\infty)\end{equation}
\begin{equation}\label{conditionsstable2}\left[ X^n,N\right]_t\stackrel{p}{\longrightarrow}0~~\forall t\,\in\,[0,\infty)~\text{and}~\forall\,N\,\in\,\mathcal{M}^{\perp}\end{equation}
\end{subequations}
are sufficient that $(X^n)$ converges $(\mathcal{F})$-stably in law to $W_{V_t}$, where $W$ is a standard Brownian motion independent of $\,\mathcal{F}$.
\end{theo}
This theorem is a simplified martingale version of the more general theorem 2--1 in \cite{jacod1}. A similar special version of the theorem is also used in \cite{fukasawa}. A comprehensive illustrative overview on Jacod's stable limit theory and further motivation and applications of this result can be found in \cite{poldilimit}. A discrete-time version of that theorem (cf.\,3--1 in \cite{jacod1}) is the following:
\begin{cor}\label{djac}Assume that $Z_t^n=\sum_{T_{n,i}\le t} X_{n,i}$ is the endpoint of a discrete martingale and the $X_{n,i}$ are $\mathcal{F}_{T_{n,i}}$-measurable square integrable random variables and $(W_t,\mathcal{F}_t)$ a Brownian motion and $\Delta T_{n,i}=T_{n,i+1}-T_{n,i}\rightarrow 0$ as $n\rightarrow\infty$. If there exists a predictable process $(v_s)_{s\ge 0}$ such that
\begin{subequations}
\begin{align}\label{djac1}\sum_{T_{n,i}\le t}\E\left[X_{n,i}^2|\mathcal{F}_{T_{n,i-1}}\right] \stackrel{p}{\longrightarrow} \int_0^tv_s^2\,ds~,\end{align}
\begin{align}\label{djac2}\forall\epsilon>0:~\sum_{T_{n,i}\le t}\E\left[X_{n,i}^2\,\1_{\{X_{n,i}>\epsilon\}}|\mathcal{F}_{T_{n,i-1}}\right]\stackrel{p}{\longrightarrow} 0~,\end{align}
\begin{align}\label{djac3}\sum_{T_{n,i}\le t}\E\left[X_{n,i}(W_{T_{n,i}}-W_{T_{n,i-1}})|\mathcal{F}_{T_{n,i-1}}\right]\stackrel{p}{\longrightarrow}0~,\end{align}
\begin{align}\label{djac4}\sum_{T_{n,i}\le t}\E\left[X_{n,i}(M_{T_{n,i}}-M_{T_{n,i-1}})|\mathcal{F}_{T_{n,i-1}}\right]\stackrel{p}{\longrightarrow}0~,\end{align}
\end{subequations}
for all bounded $\mathcal{F}_t$-martingales with $M_0=0$ and $\left[ W,M\right]\equiv 0$. Then the following stable convergence of the process $Z_t^n$ holds true:
\begin{align}Z_t^n\stackrel{st}{\rightsquigarrow}Z_t=\int_0^tv_s\,dW^{\bot}_s\end{align}
where $W^{\bot}$ is a Brownian motion defined on an orthogonal extension of the original probability space. \end{cor}
The limiting process in the foregoing Theorem \ref{jacod} is a time-changed Brownian motion. The Brownian motion is of central importance in the theory of continuous local martingales, since every continuous local martingale $M_t$ corresponds to a Dambis, Dubins-Schwarz time-changed Brownian motion $B_{\left[ M\right]_t}$. For each $(X_t^n)$ we have a representation as Dambis, Dubins-Schwarz Brownian motion $W^n_{\left[ X^n\right]_t}$ and the sequence converges weakly to a limiting Brownian motion $W_V$ by the asymptotic Knight-theorem. We refer to Theorem 7.7 in \cite{billingsley} for a proof.
The conditions \eqref{conditionsstable1} and \eqref{conditionsstable2} about the quadratic covariations converging to zero in probability ensure that the weak convergence to $W_V$ is stable.\\
For one fixed $0<T<\infty$ we have the result that $X^n_T$ converges stably in law to a centred mixed normal distribution:
\begin{equation}\label{stableclt}X^n_T\stackrel{st}{\rightsquigarrow}\mathbf{N}\left(0,V_T\right)~.\end{equation}
The independence of the limiting Brownian motion $W$ and $(V,Y)$ for any $\mathcal{F}$-measurable random variable $Y$ assures that $(W_{V_{T}},Y)$ has the same law as $(V_T Z,Y)$ with $Z\sim\mathbf{N}(0,1)$ and independent of $(V_T,Y)$.\\
Note, that in the original theorem 2--1 in \cite{jacod1} for semimartingales the same conditions as in our Theorem \ref{jacod} are imposed for the predictable quadratic (co-)variation processes that coincide with the quadratic (co-)variations for continuous semimartingales. Additionally, a condition that the drift can be neglected asymptotically is imposed. Compared to Theorem 3--1 in \cite{jacod1}, we allow for non-equidistant discrete partitions which does not harm the deduction of Theorem 3--1 from Theorem 2--1 in \cite{jacod1}.  A conditional Lindeberg-condition \eqref{djac2} and a convergence condition on the conditional variances \eqref{djac1} are analogous as in central limit theorems for triangular martingale arrays. The main difference to the stable limit theorem Corollary 3.\,1 in \cite{hall} (page 58 ff.\,) is that a certain nesting condition on the filtrations is replaced by conditions \eqref{djac3} and \eqref{djac4}. Usually the reference Brownian motion $W$ is given and ``fully generates'' the $X_{n,i}$s in the sense that \eqref{djac4} holds.\\
The theorem also extends to a multi-dimensional setting which is formulated separately in the next corollary. For this purpose let $M^{*}$ denote the transpose of a vector $M$ and the $(d\times r)$-dimensional quadratic covariation $\left[ M,N^{*}\right]_t\:=\left(\left[ M^{i},N^{j}\right]_t\right)_{ij}$ with $1\le i\le d$ and $1\le j \le r $ for a $d$-dimensional $M$ and $r$-dimensional $N$. Recall that convergence in probability of a vector is equivalent to convergence in probability for every component.
\begin{cor}\label{jacodmult}
Let $(M_t,\mathcal{F}_t)$ be a $d$-dimensional continuous local martingale and $\mathcal{M}^{\bot}$ again the set of $(\mathcal{F}_t)$-adapted bounded martingales orthogonal to $M$ (to all components). A sequence of $r$-dimensional continuous $(\mathcal{F}_t)$-adapted local martingales $(X^n)$ with
\begin{align}\label{condstablemult1}\left[ X^n,{X^n}^{*}\right]_t\stackrel{p}{\longrightarrow} V_t=\int_0^tw_s w_s^{*}\, ds~,\end{align}
where $w_s$ is a predictable $\R^r\otimes\R^r$ process, and
\begin{subequations}
\begin{align}\label{condstablemult2}\left[ X^n,M^{*}\right]_t\stackrel{p}{\longrightarrow} 0~~\forall\,t\,\in\,[0,\infty)~,\end{align}
\begin{align}\label{condstablemult3}\left[ X^n,N\right]_t\stackrel{p}{\longrightarrow}0~~\forall\,t\,\in\,[0,\infty)~and~\forall\,N\,\in\,\mathcal{M}^{\bot}~,\end{align}
\end{subequations}
converges stably in law to the process $\int_0^tw_s dW_s$, where $W$ is a $r$-dimensional standard Brownian motion independent of $\mathcal{F}$.
\end{cor}

Jacod's theorem provides a convenient stable central limit theorem for our purpose. Nesting conditions on the sequence of filtrations that are required for other stable limit theorems as in \cite{hall} and \cite{zanten} are not satisfied here.\\
Furthermore, the concept of stable convergence enables us to prove the stable weak convergence to mixed Gaussian limiting random variables under an equivalent martingale measure $\tilde \P$ after a Girsanov transformation, where the drift processes are zero. Stable convergence guarantees that the asymptotic law carries over to the case with drift under the original measure $\P$. It is in this sense commutative with measure change (cf.\,\cite{inference}). If we have the result that $Z_n\stackrel{st}{\rightsquigarrow}m+\AVAR\,\cdot\,	
\mathbf{N}(0,1)$ under $\tilde \P$ with a standard Gaussian distribution independent of $\mathcal{F}$, defined on an orthogonal extension of the original probability space and $\mathcal{F}$-measurable bounded random variables $m$ and $\AVAR$, the same convergence holds true under $\P$. Since stable convergence $Z_n\stackrel{st}{\rightsquigarrow} Z$ implies for all $f\in\mathcal{C}(\mathcal{E})$ and $\mathcal{F}$-measurable bounded random variables $X$
$$\E\left[Xf(Z_n)\right]=\tilde \E\left[(d\P/d\tilde \P)Xf(Z_n)\right]\rightarrow \tilde\E^{\prime}\left[(d\P/d\tilde \P)Xf(Z)\right]=\E^{\prime}\left[Xf(Z)\right]~,$$
by uniform integrability of $Xf(Z_n)(d\P/d\tilde \P)$ with 
$${\text{d}\P}/{\text{d} \tilde \P}=\exp{\left(-\int_0^t \gamma_s dB_s+\frac{1}{2}\int_0^t\gamma_s^2~ds\right)}~,$$
where $\sigma_s\gamma_s+\mu_s=0$.

\section{A synchronized realized covolatility estimator\label{sec:3}}
\begin{assump}\label{eff}
On a filtered probability space $\left(\Omega,\F,\left(\F_t\right),\P\right)$,  $X=(X_t)_{t\in\R^{+}}$ and $Y=(Y_t)_{t\in\R^{+}}$ are continuous semimartingales defined by the following stochastic differential equations:
\begin{align*} dX_t&=\mu_t^X\,dt+\sigma_t^X\,dB_t^X~,\\
							dY_t&=\mu_t^Y\,dt+\sigma_t^Y\,dB_t^Y~,
\end{align*}
with two $\left(\F_t\right)$--adapted standard Brownian motions $B^X$ and $B^Y$ and $\rho_t\,dt=d\left[ B^X,B^Y\right]_t$.
The drift processes $\mu_t^X$ and $\mu_t^Y$ are $\left(\F_t\right)$--adapted locally bounded stochastic processes and the spot volatilities $\sigma_t^X$ and $\sigma_t^Y$ and $\rho_t$ are assumed to be $\left(\F_t\right)$--adapted with continuous paths. We assume strictly positive volatilities and the Novikov condition $\E\left[\exp{\left((1/2)\int_0^T(\mu^{\,\cdot\,}/\sigma^{\,\cdot\,})^2_t\,dt\right)}\right]<\infty$ for $X$ and $Y$.
\end{assump}
We consider the estimation of the quadratic covariation $\left[ X,Y\right]_T$ of two continuous semimartingales, also called It\^{o} processes, $X$ and $Y$ as defined in Assumption \ref{eff} from discrete observations following non-synchronous sampling schemes.\\
We impose the following regularity assumptions on the underlying asynchronous sampling schemes:
\begin{assump}\label{grid}
The deterministic observation times $\mathcal{T}^{X,n}=\{0\le t_0^{(n)}<t_1^{(n)}<\ldots<t_n^{(n)}\le T\}$ of $X$ and $\mathcal{T}^{Y,m}=\{0\le \tau_0^{(m)}<\tau_1^{(m)}<\ldots<\tau_m^{(m)}\le T\}$ of $Y$ are assumed to be regular in the following sense:
There exists a constant $0<\alpha\le 1/3$ such that
\begin{subequations}
\begin{align}
\delta_n^X&=\sup_{i\in \{1,\ldots,n\}}{\left(\left(t_i^{(n)}-t_{i-1}^{(n)}\right),t_0^{(n)},T-t_n^{(n)}\right)}~\;\,=\mathcal{O}\left(n^{-\nicefrac{2}{3}-\alpha}\right)~,\\
\delta_m^Y&=\sup_{j\in \{1,\ldots,m\}}{\left(\left(\tau_j^{(m)}-\tau_{j-1}^{(m)}\right),\tau_0^{(m)},T-\tau_{m}^{(m)}\right)}=\mathcal{O}\left(m^{-\nicefrac{2}{3}-\alpha}\right)~.
\end{align}
\end{subequations}
We consider asymptotics where the number of observations of $X$ and $Y$ are assumed to be of the same asymptotic order $n=\mathcal{O}(m)$ and $m=\mathcal{O}(n)$ and express that shortly by $n\sim m$.
\end{assump}
For synchronous data $n=m$ and $t_i^{(n)}=\tau_i^{(n)}$ for all $i\in\{0,\ldots,n\}$ holds. In the non-synchronous case the number of observations $(n+1)$ of $X$ and $(m+1)$ of $Y$ may differ and the sets of observation times $\mathcal{T}^{X,n}$ also contain times $t_i^{(n)}\notin\mathcal{T}^{Y,m}$ and $\tau_j^{(n)}\notin\mathcal{T}^{X,n}$. We work within the general model where also synchronous observation times can take place and hence $\mathcal{T}^{Y,m}$ and $\mathcal{T}^{X,n}$ are not assumed to be disjoint. In the following, we omit the superscripts $(n)$ and $(m)$ for observation times to increase the readability.\\
Although the sequences of observation times are modeled deterministically, we remark that the case of random sampling times that are independent of the observed processes is included in that analysis regarding the conditional law given the observation times. \\
We use the short notation $\d X_{t_i},i=1,\ldots,n$ from now on for increments $X_{t_i}-X_{t_{i-1}}$ and analogously for $Y$. In \cite{hy} the consistency of the estimator
$$\widehat{\left[  X, Y\right]}_T^{(HY)}=\sum_{i=1}^{n}\sum_{j=1}^{m}\d  X_{t_i}\d  Y_{\tau_j}\1_{[\min{(t_{i},\tau_{j})}>\max{(t_{i-1},\tau_{j-1})}]}~,$$
is proved, where the product terms include all increments of the processes with overlapping observation time intervals, for a similar model of discretely observed It\^{o} diffusions with deterministic correlation, drift and volatility functions. Consistency directly carries over to our setting including random correlation, drift and volatility processes. The estimator is also in our setting, furthermore, unbiased if drift terms are zero and else asymptotically unbiased. In \cite{hy2} it has further been shown that on stronger regularity assumptions on the observation schemes this Hayashi-Yoshida estimator is asymptotically distributed according to a Gaussian law.\\
For a general strategy leading to a synchronization mechanism that keeps to the Hayashi-Yoshida approach and its valuable properties, we focus on an alternative useful method to handle the asynchronicity of the data. It has been introduced in \cite{palandri}, where it was called pseudo-aggregation. The method translates the Hayashi-Yoshida estimator into an iterative algorithm that allows to rewrite the estimator without indicator functions. This can be done by aggregation of addends for which partial sums are telescoping. A first simple rewriting of the Hayashi-Yoshida estimator is obtained by taking the sum of the products of all increments of $X$ with the telescoping sums of aggregated observed increments of $Y$ for that observation time instants overlap with the according observation time instant of $X$ (or in the symmetric way):
\begin{align*}\widehat{\left[  X, Y\right]}_T^{(HY)}&=\sum_{i=1}^{n}\d  X_{t_i}\left(\sum_{j\in\{1,\ldots,m\}}\d  Y_{\tau_j}\1_{[\min{(t_{i},\tau_{j})}>\max{(t_{i-1},\tau_{j-1})}]}\right)\\
&=\sum_{j=1}^{m}\d  Y_{\tau_j}\left(\sum_{i\in\{1,\ldots,n\}}\d  X_{t_i}\1_{[\min{(t_{i},\tau_{j})}>\max{(t_{i-1},\tau_{j-1})}]}\right)~.
\end{align*}

\renewcommand{\figurename}{Algorithm}
\begin{figure}[!ht]\fbox{
\begin{minipage}[c]{\textwidth} 
\renewcommand{\baselinestretch}{.75}\normalsize
first step:
\begin{itemize}\item for $t_{0}<\tau_{0}$ and $\mu_0=\min{(w \in\{1,\ldots,n\}|\tau_{0}\le t_w)}$:
$$\mathcal{H}^0=\{t_{0},\ldots,t_{\mu_0}\}~~\mbox{and}~~\mathcal{G}^0=\{\tau_{0}\}$$
$$q_{1}= \begin{cases}\mu_0+1~~~\mbox{if}~~\tau_{0}=t_{\mu_0}\\ \mu_0 ~~~~~~~~~~\mbox{if}~~\tau_{0}<t_{\mu_0}\end{cases}~~\mbox{and}~~ r_{1}=1$$
\item for $t_{0}=\tau_{0}$:
$$\mathcal{H}^0=\{t_{0}\}~~\mbox{and}~~\mathcal{G}^0=\{\tau_{0}\}$$
$$q_{1}=1~~~\mbox{and}~~~r_1=1$$
\item for $t_{0}>\tau_{0}$ and $w_0=\min{(l \in \{1,\ldots,m\}|t_{0}\le \tau_l)}$:
$$\mathcal{H}^0=\{t_{0}\}~~\mbox{and}~~\mathcal{G}^0=\{\tau_{0},\ldots,\tau_{w_0}\}$$
$$q_{1}=1~~\mbox{and}~~r_1=\begin{cases}w_0+1~~~\mbox{if}~~t_{0}=\tau_{w_0}\\ w_0 ~~~~~~~~~~\mbox{if}~~t_{0}<\tau_{w_0}\end{cases}$$
\end{itemize}
$i$th step (given $\mathcal{H}^{i-1}$ and $\mathcal{G}^{i-1}$):
\begin{itemize}\item for $t_{q_i}<\tau_{r_i}$ and $\mu_i=\min{(w \in\{q_i+1,\ldots,n\}|\tau_{r_i}\le t_w)}$:
$$\mathcal{H}^i=\{t_{q_i},\ldots,t_{\mu_i}\}~~\mbox{and}~~\mathcal{G}^i=\{\tau_{r_i}\}$$
$$q_{i}\dashrightarrow\begin{cases}q_{i+1}=\mu_i+1~~~\mbox{if}~~\tau_{r_i}=t_{\mu_i}\\ q_{i+1}=\mu_i ~~~~~~~~~~\mbox{if}~~\tau_{r_i}<t_{\mu_i}\end{cases}~~\mbox{and}~~ r_{i}\dashrightarrow r_{i+1}=r_i+1$$
\item for $t_{q_i}=\tau_{r_i}$:
$$\mathcal{H}^i=\{t_{q_i}\}~~\mbox{and}~~\mathcal{G}^i=\{\tau_{r_i}\}$$
$$q_{i}\dashrightarrow q_{i+1}=q_i+1~~~\mbox{and}~~~r_i\dashrightarrow r_{i+1}=r_i+1$$
\item for $t_{q_i}>\tau_{r_i}$ and $w_i=\min{(l \in \{r_i+1,\ldots,m\}|t_{q_i}\le \tau_l)}$:
$$\mathcal{H}^i=\{t_{q_i}\}~~\mbox{and}~~\mathcal{G}^i=\{\tau_{r_i},\ldots,\tau_{w_i}\}$$
$$q_{i}\dashrightarrow q_{i+1}=q_i+1~~\mbox{and}~~r_i\dashrightarrow\begin{cases}r_{i+1}=w_i+1~~~\mbox{if}~~t_{q_i}=\tau_{w_i}\\ r_{i+1}=w_i ~~~~~~~~~~\mbox{if}~~t_{q_i}<\tau_{w_i}\end{cases}$$
\end{itemize}
\end{minipage}}
\caption{\label{A}Iterative algorithm for construction of the joint grid from asynchronous data.}
\end{figure}
\renewcommand{\figurename}{Figure}
\addtocounter{figure}{-1} 
\renewcommand{\baselinestretch}{1}\normalsize\noindent
Defining the next-tick interpolation $t_{i,+}\:=\min_{0\le j\le m}{\left(\tau_j|\tau_j\ge t_i\right)}$ and the previous-tick interpolation $t_{i,-}\:=\max_{0\le j\le m}{\left(\tau_j|\tau_j\le t_i\right)}$, the last expression can be illustrated
$$\widehat{\left[  X, Y\right]}_T^{(HY)}=\sum_{i=1}^n\d X_{t_i}\left(Y_{t_{i,+}}-Y_{t_{i-1,-}}\right)~.$$
The algorithm which we will use is a more enhanced method to aggregate the data in an adequate way. For this purpose $(N+1)$ sets $\mathcal{H}^i$ and $\mathcal{G}^i$ are constructed, where $N<\min{(n,m)}$, each set including one or more than one observation time of $X$ and $Y$, respectively. This method to construct a joint grid for the observations of the two processes is described by Algorithm \ref{A}.\\
The Algorithm \ref{A} that we have first presented in \cite{bibinger} stops after $(N+1)$ steps when the last observation time is reached. We pass over from the original observations to the sums of observed increments $X^{\mathcal{H}^{i}}$ over sets $\mathcal{H}^i$ and $Y^{\mathcal{G}^{i}}$ over sets $\mathcal{G}^i$, respectively. The observations are grouped together so that the resulting realized covolatility estimator 
$$\sum_{i=0}^{N} X^{\mathcal{H}^{i}} Y^{\mathcal{G}^{i}}=\sum_{i=1}^{n}\sum_{j=1}^{m}\d X_{t_i}\d  Y_{\tau_j}\1_{[\min{(t_{i},\tau_{j})}>\max{(t_{i-1},\tau_{j-1})}]}$$
calculated from the `synchronized' observations 
$$X^{\mathcal{H}^{i}}=\sum_{{t_j}\in \mathcal{H}^i}\d  X_{t_j}~,~~Y^{\mathcal{G}^{i}}=\sum_{{\tau_j}\in \mathcal{G}^i}\d  Y_{\tau_j}~,~~i \in \{0,\ldots,N\}~.$$
for the integrated covolatility will coincide with the one by \cite{hy} stated above.
We use a different illustration of this estimator compared to \cite{palandri} making use of telescoping sums.\\
With the denotation expressions from Algorithm \ref{A} 
\begin{align*}&\mu_i=\max{(k|t_k \in \mathcal{H}^i)},&w_i=\max{(k|\tau_k \in \mathcal{G}^i)}&~~\mbox{and}\\ &q_i=\min{(k|t_k \in \mathcal{H}^i)},&r_i=\min{(k|\tau_k \in \mathcal{G}^i)}\,&~~,i\in\{0,\ldots,N\}\end{align*} 
and for the purpose of a simpler notation
\begin{align*}& X_{g_i}= X_{t_{\mu_i}},&Y_{ \gamma_i}= Y_{\tau_{w_i}}~~~~\,&,i\in\{0,\ldots,N\}~\mbox{and}\\
& X_{l_i}\,=X_{t_{q_i-1}},& Y_{\lambda_i}= Y_{\tau_{r_i-1}}~~&,i\in\{1,\ldots,N\}\end{align*}
with $l_0\:=t_0,\,\lambda_0\:=\tau_0$, $ X^{\mathcal{H}^{i}}$ and $ Y^{\mathcal{G}{i}}$ can be written as telescoping sums 
$X^{\mathcal{H}^{i}}=\left( X_{g_i}- X_{l_i}\right)$,
 $Y^{\mathcal{G}{i}}=\left( Y_{\gamma_i}-Y_{\lambda_i}\right)~.$
This leads to
\begin{align}\label{HY}\widehat{\left[  X, Y\right]}_T^{(HY)}=\sum_{i=1}^N\left( X_{g_i}- X_{l_i}\right)\left( Y_{\gamma_i}- Y_{\lambda_i}\right)~,\end{align}
where summation starts with $i=0$ or $i=1$ since the addend for $i=0$ is always zero.
Although we use this specific new illustration throughout this article, we will call this realized covolatility of our synchronized observations also  Hayashi-Yoshida estimator in the following.
In this notation $g_i$ denotes the greatest and $l_i$ the last observation time before the least element of the set $\mathcal{H}^i$ and analogously $\gamma_i$ and $\lambda_i$ of $\mathcal{G}^i$.\\ \\
%\rule{\textwidth}{.5pt}\\
\textbf{Example}\\
An illustration of the application of Algorithm \ref{A} to observations is given in Figure \ref{example}. In this example, we have
$\mathcal{H}^0=\{t_0\},\mathcal{G}^0=\{\tau_0\}, 
\mathcal{H}^1=\{t_1,t_2,t_3\},\mathcal{G}^1=\{\tau_1\},
\mathcal{H}^2=\{t_3\},\mathcal{G}^2=\{\tau_2,\tau_3\},
\mathcal{H}^3=\{t_4,t_5,t_6\},\mathcal{G}^3=\{\tau_4\},
\mathcal{H}^4=\{t_6,t_7\},\mathcal{G}^4=\{\tau_5\},
\mathcal{H}^5=\{t_7,t_8\},\mathcal{G}^5=\{\tau_6\},
\mathcal{H}^6=\{t_8\},\mathcal{G}^6=\{\tau_7,\tau_8\},
\mathcal{H}^7=\{t_9\},\mathcal{G}^7=\{\tau_8,\tau_9\},
\mathcal{H}^8=\{t_{10}\},\mathcal{G}^8=\{\tau_9,\tau_{10}\}$~.\\

\begin{figure}[!t]
\begin{center}
\fbox{\includegraphics[width=14.4cm]{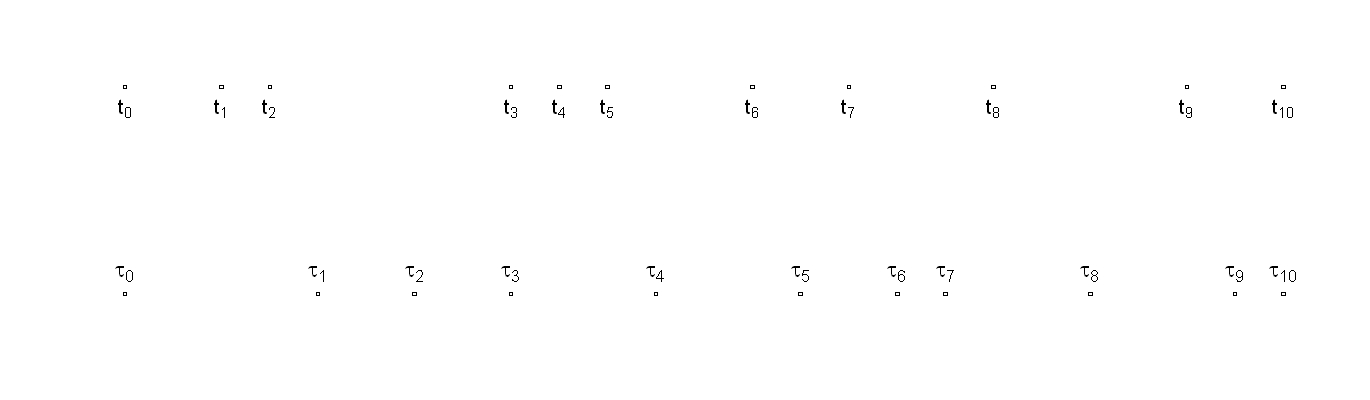}}\end{center}
\caption{\label{example}Example for synchronization using Algorithm \ref{A}.}
\end{figure}\noindent
The example highlights the important features of the synchronization procedure. The sets $\mathcal{H}^i$ and $\mathcal{G}^i$ are in general not disjoint and the maxima of consecutive sets can be the same time points. The minimum of a successive set can as well equal the maximum of the prevenient. Contrarily, consecutive minima are not equal. For further examples we refer to \cite{palandri}. Of course the example is just for illustration and the number of observations is much smaller than in practice.
The synchronization of $n+1=11$ and $m+1=11$ observations leads to $N+1=9$ synchronized observations in this example.\\ \\
%\rule{\textwidth}{.5pt}\\ \\
The fact that we obtain $(N+1)<\min{(n,m)}+1$ synchronized observations indicates heuristically that the efficiency of such techniques of covariance estimation mainly depends on the number of observations available for the less liquid process which is observed at a lower frequency. By Assumption \ref{grid} we restrict us to the case that $n$ and $m$ are of the same order. Thus for the suprema of times between two observations 
$$\delta_n^X=\mathcal{O}\left(N^{-\nicefrac{2}{3}-\alpha}\right)~~\text{and}~~\delta_n^Y=\mathcal{O}\left(N^{-\nicefrac{2}{3}-\alpha}\right)$$
holds with a constant $0<\alpha\le 1/3$.\\
In the next section, we show that on Assumption \ref{eff} and \ref{grid} the estimator \eqref{HY} is $\sqrt{N}$-consistent and, on further assumptions on the asymptotic behavior of the asynchronous sampling schemes, asymptotically normally distributed. Using standard interpolation methods such an estimator cannot be obtained.\\
Another recent approach to deal with non-synchronous discrete observations in a general setting including market microstructure noise has been proposed by \cite{bn1}. This method is also related to our approach. The so-called refresh times are the cumulative sums of waiting times until both processes are observed. Assume that in the $i$th step of Algorithm \ref{A} $t_{q_i}<\tau_{r_i}$ holds. Then the next observation times of $X$ are grouped together ending with the first observation time $t_{\mu_i-1}<\tau_{r_i}\le t_{\mu_i}$ greater or equal than $\tau_{r_i}$. Then we start the next comparison step and compare this last observation time grouped to the set $\mathcal{H}^{i}$ to $\tau_{r_i+1}$, except for the case where two synchronous observations appeared, where we compare the two following times. Since in the completely asynchronous case at the refresh times only one of the two processes is observed, the refresh time method used in \cite{bn1} includes a previous-tick interpolation for the unobserved process at the refresh times. Refresh times provide the `closest synchronous approximation' to the asynchronous sampling schemes that we define in Proposition \ref{partition} below. The number of refresh times which are denoted in this work by $T_i, i=0,\ldots,N$, equals the number of sets constructed by pseudo-aggregation. In a setting that also takes microstructure noise into account, a consistent estimator requires smoothing techniques to reduce the noise perturbation and the optimal convergence rate is slower (cf.\,\cite{bibinger}). The previous-tick interpolation, however, causes a negative bias due to asynchronicity when calculating the simple realized covolatility estimator based on the refresh time and previous-tick approach and it does not equal the estimator of Hayashi-Yoshida. The reason for this bias is that, due to the previous-tick interpolation, products of increments with overlapping observation time instants fall out of the realized covolatility. The pseudo-aggregation Algorithm \ref{A} used in this work corresponds to the refresh time method when replacing the previous-tick interpolation by a next-tick interpolation for the right end points of refresh time instants. Then, the resulting realized covolatility of `synchronized observations' 
\begin{align}\notag\widehat{\left[  X, Y\right]}_T^{(HY)}&=\sum_{i=1}^N\left( X_{g_i}- X_{l_i}\right)\left( Y_{\gamma_i}- Y_{\lambda_i}\right)\\
&=\sum_{i=1}^N\left( X_{T_{i,+}^X}- X_{T_{i-1,-}^X}\right)\left( Y_{T_{i,+}^Y}- Y_{T_{i-1,-}^Y}\right)\end{align}
coincides with the Hayashi-Yoshida estimator and has no bias due to asynchronicity. As figured out in the simulation study of \cite{bibinger} the asymptotically vanishing influence of the bias due to pure previous-tick interpolation also shows up in the setting with noise for finite sample sizes and mild noise variances for that combined estimators are constructed in \cite{bn1} and \cite{bibinger}, among others.

\begin{figure}[t]
\begin{center}\hspace*{-.75cm}
\fbox{\includegraphics[width=14.4cm]{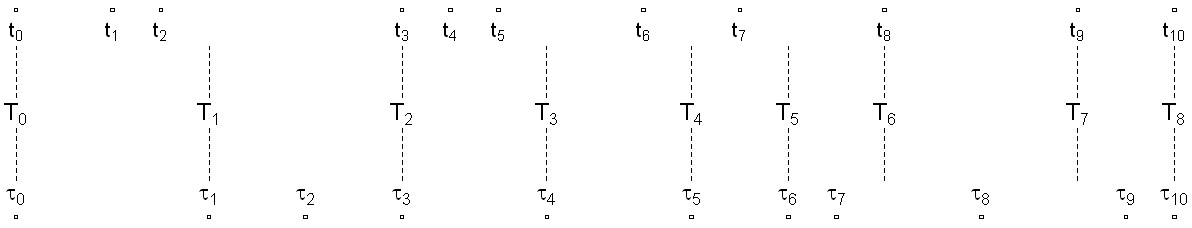}}\end{center}
\caption{\label{example2}Example for synchronization using Algorithm \ref{A} including refresh times.}
\end{figure}\noindent
Figure \ref{example2} visualizes refresh times $T_i,\,i=0,\ldots,8$ for our above given example. For this example the realized covolatility calculated with refresh time previous-tick interpolated values equals
\begin{align*}(X_{t_2}-X_{t_0})(Y_{\tau_1}-Y_{\tau_0})+(X_{t_3}-X_{t_2})(Y_{\tau_3}-Y_{\tau_1})+(X_{t_5}-X_{t_3})(Y_{\tau_4}-Y_{\tau_3})+\\(X_{t_6}-X_{t_5})(Y_{\tau_5}-Y_{\tau_4})+(X_{t_7}-X_{t_6})(Y_{\tau_6}-Y_{\tau_5})+(X_{t_8}-X_{t_7})(Y_{\tau_7}-Y_{\tau_6})+\\(X_{t_9}-X_{t_8})(Y_{\tau_8}-Y_{\tau_7})+(X_{t_{10}}-X_{t_9})(Y_{\tau_{10}}-Y_{\tau_8})\end{align*}
and is biased downwards due to non-synchronicity, whereas \eqref{HY} yields
\sectionmark{Asymptotics of the HY-estimator}
\begin{align*}(X_{t_3}-X_{t_0})(Y_{\tau_1}-Y_{\tau_0})+(X_{t_3}-X_{t_2})(Y_{\tau_3}-Y_{\tau_1})+(X_{t_6}-X_{t_3})(Y_{\tau_4}-Y_{\tau_3})+\\(X_{t_7}-X_{t_5})(Y_{\tau_5}-Y_{\tau_4})+(X_{t_8}-X_{t_6})(Y_{\tau_6}-Y_{\tau_5})+(X_{t_8}-X_{t_7})(Y_{\tau_8}-Y_{\tau_6})+\\(X_{t_9}-X_{t_8})(Y_{\tau_9}-Y_{\tau_7})+(X_{t_{10}}-X_{t_9})(Y_{\tau_{10}}-Y_{\tau_8})~,\end{align*}
which is an unbiased estimator for observations of processes according to Assumption \ref{eff}, when drift terms are assumed to be zero.

\section{Asymptotic distribution theory\label{sec:4}}
\begin{figure}[t]
\begin{center}
\includegraphics[width=10cm]{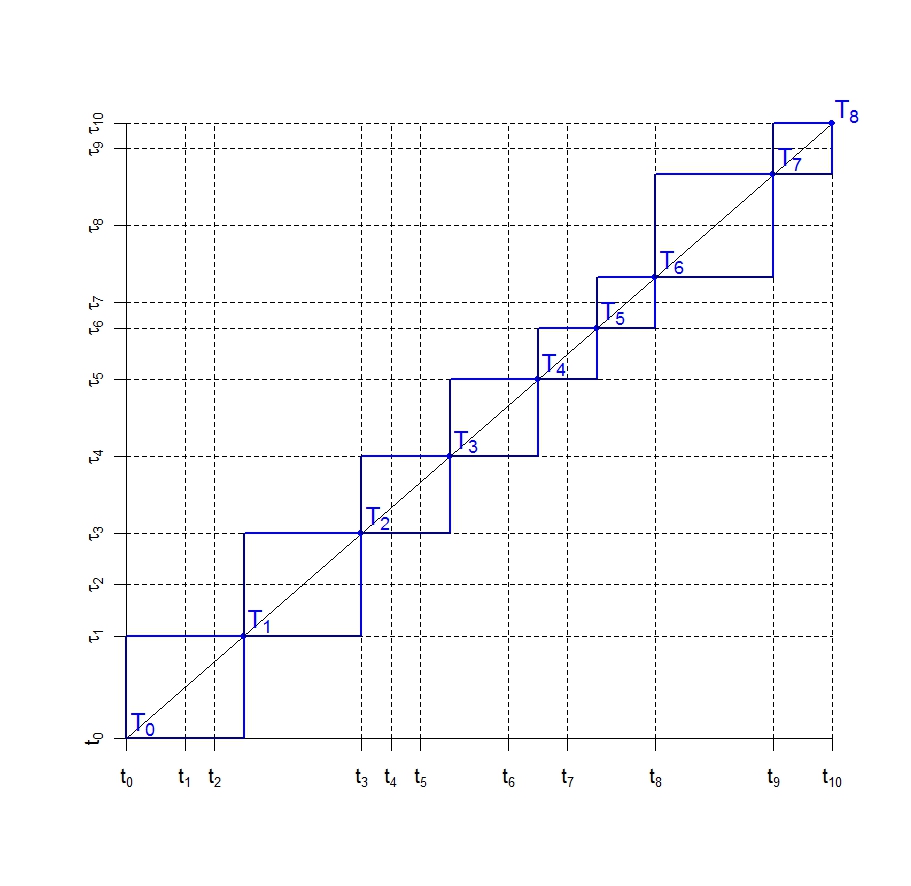}\end{center}
\vspace*{-1cm}
\caption{\label{quader1}Illustration of the synchronous approximation for our example.}
\end{figure}

\begin{figure}[ht]
\begin{center}
\includegraphics[width=10cm]{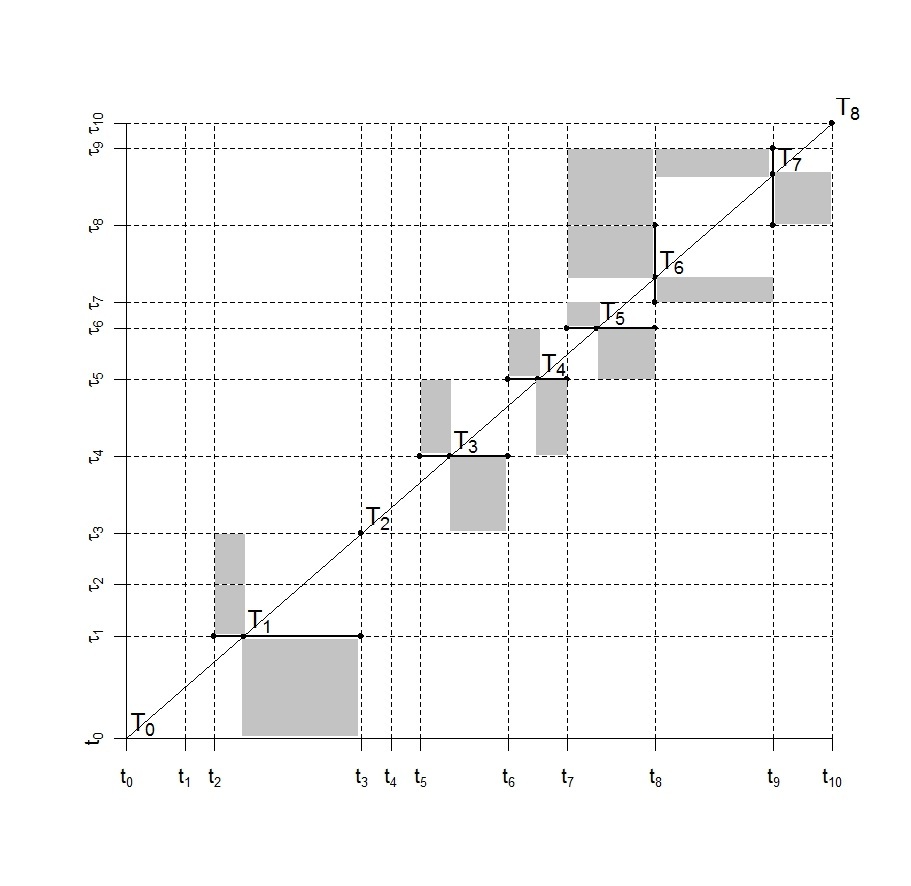}\end{center}
\vspace*{-1cm}
\caption{\label{quader4}Illustration of the next- and previous-tick interpolated values and the error due to non-synchronicity for our example.}
\end{figure}

\begin{figure}[t]
\framebox{
\includegraphics[width=7.285cm]{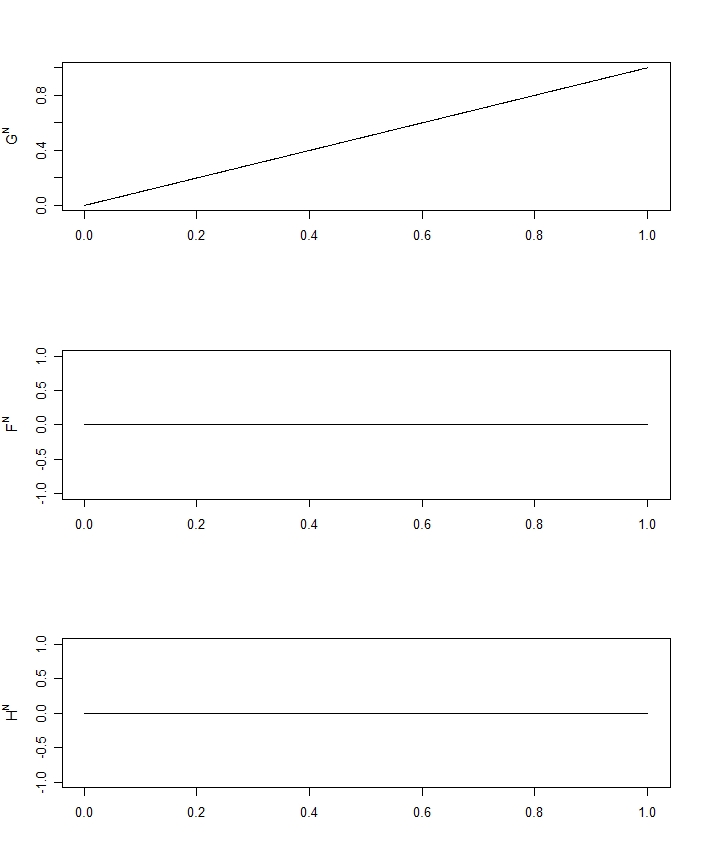}
\includegraphics[width=7.285cm]{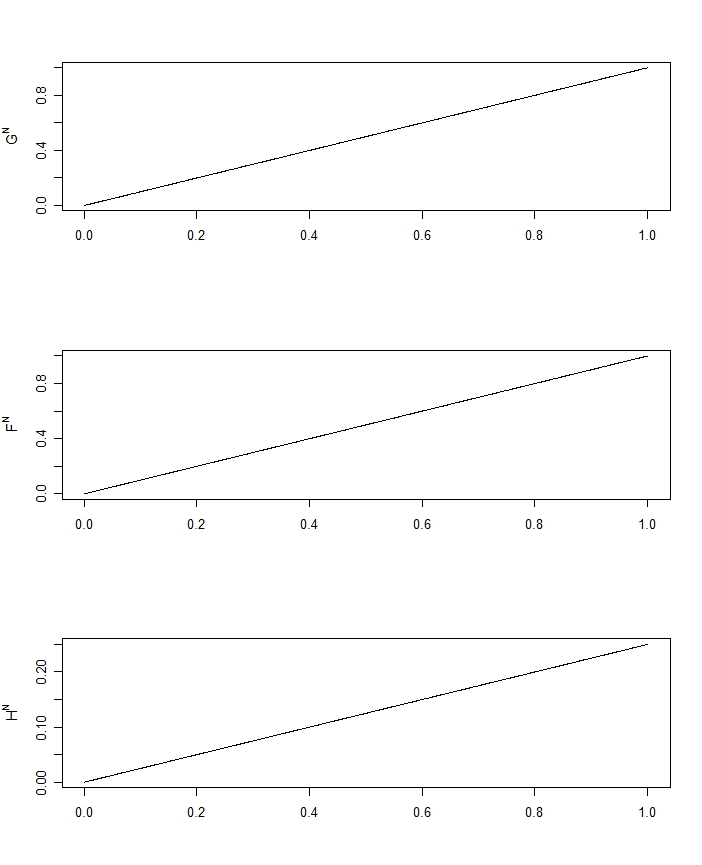}
}
\caption{\label{figqvt}Quadratic (co-)variations of time for synchronous equidistant (left) and intermeshed (right) sampling.} 
\end{figure}

In this section the elements for an analysis of the asymptotic properties of the estimator \eqref{HY} are developed where the emphasis is on the asymptotic distribution of the estimator.\\
The following technical Proposition constitutes the theoretical justification that the refresh times $T_i^{(n)},\,1\le i\le N$ introduced in the foregoing section can serve as a convenient basis to decompose the overall estimation error of the synchronized realized covolatility \eqref{HY}. For every $N$ these times induce a partition of the time horizon $[0,T]$ that we call the closest synchronous approximation.
\begin{prop}\label{proppartition}
If we define $T_i^{(N)}\:=\min{\left(g_i,\gamma_i\right)},~i=0,\ldots,N$, the set $\mathcal{T}^{syn,N}=\{T_0^{(N)},\ldots,T_N^{(N)}\}$ induces a partition of the time span \,$[0,T]$ in the sense that $\dunion_i [T_i^{(N)},T_{i+1}^{(N)})=[T_0^{(N)},T-T_N^{(N)})$.\\
The following equality holds true:
\begin{equation}\label{partition} T_i^{(N)}=\min{\left(g_i^{(N)},\gamma_i^{(N)}\right)}=\max{\left(l_{i+1}^{(N)},\lambda_{i+1}^{(N)}\right)},~i=1,\ldots,N-1\end{equation}
and on Assumption \ref{grid} $\delta_N\:=\sup_{i\in\{1,\ldots,N\}}{\left(T_i^{(N)}-T_{i-1}^{(N)}\right)}=\mathcal{O}\left(N^{-\nicefrac{2}{3}-\alpha}\right)$ holds.
\end{prop}
In the following we frequently leave out superscripts indicating dependence on $N$ to guarantee clarity and increase the readability.
\begin{proof} Assume without loss of generality $g_i\le\gamma_i$ for an arbitrarily fixed $i\in\{1,\ldots,N-1\}$. Taking Algorithm \ref{A} into account, we proof that \eqref{partition} holds true.\\
If $g_i<\gamma_i$, then the observation times $\gamma_i$ and $g_{i,+}\:=\min{\left(t_k\in \mathcal{T}^X|t_k>g_i\right)}$ are compared in the $(i+1)$th step of the synchronization Algorithm \ref{A} and $g_{i,+}=\min{\left(t_k\in \mathcal{T}^X|t_k\in\H^{i+1}\right)}$ holds true. Thus, $g_i=l_{i+1}$ and \eqref{partition} holds true. We remark that in this case $\gamma_i\in \G^{i+1}$ and thus $\gamma_i>\lambda_{i+1}=\gamma_{i,-}\:=\max{\left(\tau_k\in\mathcal{T}^Y|\tau_k<\gamma_i\right)}\ge \gamma_{i-1}$.\\
If $g_i=\gamma_i$, then the observation times $g_{i,+}$ and $\gamma_{i,+}$ are compared in the $(i+1)$th step of Algorithm \ref{A} and $l_{i+1}=\lambda_{i+1}=g_i=\gamma_i$ what implies \eqref{partition}.\\ Equation \eqref{partition} does not hold true for $i=0,N$ and $T_0=t_0\wedge\tau_0$ because we have set $l_0=t_0$ and $\lambda_0=\tau_0$.\\
Although consecutive maxima $g_i$ of the sets $\H^{i}$ and $\gamma_i$ of the sets $\G^{i}$, respectively, can be equal, $T_i>T_{i-1}$ holds for all $i\in \{1,\ldots,N\}$ because $g_{i+1}=g_i$ implies that $\gamma_{i+1}>\gamma_i$ and $\gamma_{i+1}=\gamma_i$ implies that $g_{i+1}>g_i$. Hence, the set $\mathcal{T}^{syn}$ induces a partition of the time span $[0,T]$.
\end{proof}
The times $T_i,i=0,\ldots,N$ defined through \eqref{partition} equal the refresh times from \cite{bn1} as has been mentioned in the last section. We use Proposition \ref{proppartition} to split the error of the estimator \eqref{HY} for the integrated covolatility $\left[ X,Y\right]_T$ in two asymptotically uncorrelated parts. The error of the estimator \eqref{HY} can be written
\begin{align*}\sum_{i=1}^N\left(X_{g_i}-X_{l_i}\right)\left(Y_{\gamma_i}-Y_{\lambda_i}\right)-\int_0^T\rho_t\sigma_t^X\sigma_t^Y\,dt=D_T^N+A_T^N\end{align*}
where
\begin{align}\label{D}D_T^N\:=\sum_{i=1}^N\left(\left(X_{T_i}-X_{T_{i-1}}\right)\left(Y_{T_i}-Y_{T_{i-1}}\right)-\int_{T_{i-1}}^{T_i}\rho_t\sigma_t^X\sigma_t^Y\,dt\right)\\
\notag~~~-\int_0^{t_0\wedge \tau_0}\rho_t\sigma_t^X\sigma_t^Y\,dt-\int_{t_n\wedge\tau_m}^T\rho_t\sigma_t^X\sigma_t^Y\,dt\end{align}
is a synchronous-type discretization error of the realized covolatility estimator evaluated with synchronous observations at the times $T_i,i=0,\ldots,N$, which is the closest synchronous approximation to the asynchronous sampling scheme, and
\begin{align}\label{As}\notag A_T^N=\sum_{i=1}^N\left(Y_{\gamma_i}-Y_{\lambda_i}\right)\left(X_{g_i}-X_{T_i}\right)\1_{\{T_i=\gamma_i\}}+\left(Y_{T_{i}}-Y_{T_{i-1}}\right)\left(X_{T_{i-1}}-X_{l_i}\right)\1_{\{T_{i-1}=\lambda_i\}}\\
~~+\sum_{i=1}^N\left(X_{T_i}-X_{l_{i}}\right)\left(Y_{\gamma_i}-Y_{T_i}\right)\1_{\{T_i=g_i\}}+\left(X_{T_{i}}-X_{T_{i-1}}\right)\left(Y_{T_{i-1}}-Y_{\lambda_i}\right)\1_{\{T_{i-1}=l_i\}}\end{align}
is the remaining additional error due to the lack of synchronicity. 
When we write the increments involved in the estimator \eqref{HY} in the way
$$\left(X_{g_j}-X_{l_j}\right)=X_j^++X_j^S+X_j^-~,~\left(Y_{\gamma_j}-Y_{\lambda_j}\right)=Y_j^++Y_j^S+Y_j^-~,$$
where $X_j^+=X_{g_j}-X_{T_j}$ denotes the next-tick interpolation error at right-end points, $X_j^-=X_{T_{j-1}}-X_{l_j}$ the previous-tick interpolation error at left-end points, $X_j^S=X_{T_j}-X_{T_{j-1}}~,~j=1,\ldots,N$ the increment over the time instant of the closest synchronous approximation and analogously for $Y$, $D_T^N$ and $A_T^N$ can be expressed:
\begin{align*}D_T^N&=\sum_{i=1}^N X_i^SY_i^S~,\\
A_T^N&=\sum_{i=1}^N\left(X_i^+ (Y_i^S+Y_i^-)+Y_i^+(X_i^S+X_i^-)+X_i^-Y_i^S+Y_i^-X_i^S\right)~.\end{align*}
$D_T^N$ is an usual synchronous-type realized covolatility but incorporates an idealized sampling design at the times of the closest synchronous approximation for which we do not have observations in an asynchronous setting. Nevertheless, this idealized approximation turns out to be helpful for our further analysis. The error due to non-synchronicity $A_T^N$ hinges on the interpolations that have to be carried out since we do not observe $X$ and $Y$ at the times $T_i^N,1\le i\le N$.
The term is asymptotically centred since only products of increments over disjoint time instants remain whereas $D_T^N$ is an unbiased estimator for $\left[X,Y\right]_T$. Since either $X$ or $Y$ is observed at a certain $T_i, 1\le i\le N$, one of each interpolation errors in the illustration above equals zero.
\begin{prop}
The Brownian parts of $A_T^N$ and $D_T^N$ are uncorrelated. This means, that if we assume the drift terms to be identically zero in Assumption \ref{eff}, $A_T^N$ and $D_T^N$ are uncorrelated. If the drift terms are non-zero, $A_T^N$ and $D_T^N$ are asymptotically uncorrelated.
\end{prop}
\begin{proof}
$A_T^N$ and $D_T^N$ are both centred. If Assumption \ref{eff} holds with $\mu_t^X\equiv \mu_t^Y\equiv 0$, the expectation of the product of $A_T^N$ and $D_T^N$ is zero, since the previous- and next-tick interpolated increments in \eqref{As} are centred and uncorrelated to the other three factors in each addend of the inner sums. \\
If we allow for non-zero drift terms, Assumption \ref{eff} and Assumption \ref{grid} ensure that the increments over time intervals due to the drift induce terms at most of order $\delta_N$ in probability by products of drift terms and at most of order $\delta_N^{\nicefrac{1}{2}}$ in probability by products of drift and Brownian increments in the overall correlation.
\end{proof}
In Figure \ref{quader1} the observation times $\tau_j,j=0,\ldots,11$ of $Y$ for our Example \ref{example} from the last section are plotted against the observation times $t_i,i=0,\ldots,11$ of $X$. The dashed lines intersect for synchronous observation times $t_0=\tau_0,t_3=\tau_3$ and $t_{10}=\tau_{10}$ on the diagonal of the square in Figure \ref{quader1}. A similar visualization of the realized covolatility estimator for synchronous and equidistant data would yield coextensive squares around the diagonal, over which multiplied increments are summed up. Refresh times are (in general) not equidistant but provide a synchronous realized covolatility estimator as an approximation.
The Hayashi-Yoshida estimator \eqref{HY} is the sum of products of increments with overlapping observation time instants. The relation to the synchronous approximation $D_T^N$ is that we have next-tick interpolations and previous-tick interpolations to the times $T_i,i=0,\ldots,8$ and take increments from previous-tick interpolated values to next-tick interpolated values. The time instants of $D_T^N$ are visualized for our example in Figure \ref{quader1}. The previous- and next-tick interpolations are illustrated in Figure \ref{quader4}. The products of time instants leading to the error $A_T^N$ are illustrated in the same picture by the grey rectangles. As can be seen for the example in Figure \ref{quader4}, $A_T^N$ is the sum of the errors by the $i$th next-tick interpolation multiplied with the increments of the other process over $[\min{(l_i,\lambda_i)},T_i]$ and the sum of the errors of the $i$th previous-tick interpolation multiplied with the increments of the other process over $[T_{i-1},T_i]$.
The sum of the increments over the squares in Figure \ref{quader1}, $D_T^{8}$ for our example, and the grey rectangles in Figure \ref{quader4}, $A_T^{8}$ for our example, is the Hayashi-Yoshida estimator evaluated at the end of the last section.
\noindent
\begin{defi}[\textbf{quadratic (co-)variations of time}]\label{qvt}
For any $N\in\N$ let $T_i^{(N)},~i=0,\ldots,N$ be the times from the partition of $[0,T]$ defined in \eqref{partition} above and $g_i^{(N)},\gamma_i^{(N)},l_i^{(N)},\lambda_i^{(N)}$ the corresponding observation times designated by Algorithm \ref{A} from the estimator \eqref{HY}. $T/N$ is the mean of the time instants $\d T_i^{(N)}=T_i^{(N)}-T_{i-1}^{(N)},~i=1,\ldots,N$. Define the following sequences of functions
\begin{subequations}
\begin{align}\label{qvtg} G^N(t)=\frac{N}{T}\sum_{T_i^{(N)}\le t}\left(\d T_i^{(N)}\right)^2~,\end{align}
\begin{align}\label{qvtf}\notag F^N(t)=\frac{N}{T}\sum_{T_{i+1}^{(N)}\le t}(T_i^{(N)}-\lambda_i^{(N)})(g_i^{(N)}-T_i^{(N)})+\left(T_i^{(N)}-l_{i}^{(N)}\right)\left(\gamma_{i}^{(N)}-T_i^{(N)}\right)\\+\d T_{i+1}^{(N)}\left(T_i^{(N)}-l_{i+1}^{(N)}\right)+\d T_{i+1}^{(N)}\left(T_i^{(N)}-\lambda_{i+1}^{(N)}\right)~,\end{align}
\begin{align}\label{qvth} H^N(t)=\frac{N}{T}\sum_{T_{i+1}^{(N)}\le t}\left(T_i^{(N)}-l_{i+1}^{(N)}\right)\left(g_i^{(N)}-T_i^{(N)}\right)+\left(T_i^{(N)}-\lambda_{i+1}^{(N)}\right)\left(\gamma_i^{(N)}-T_i^{(N)}\right)~,\end{align}
\end{subequations}
for $t\in[0,T]$ that we call sequences of quadratic (co-)variations of times. 
\end{defi}
A stable central limit theorem for the estimation error is deduced on the assumption that the sequences defined by \eqref{qvtg}, \eqref{qvtf} and \eqref{qvth} converge pointwise and the sequences of difference quotients uniformly:
\begin{assump}[\textbf{asymptotic quadratic (co-)variation of times}]\label{aqvt}
Assume that for the sequences of sampling schemes and the times $T_i^{(N)},g_i^{(N)},\gamma_i^{(N)},l_i^{(N)},\lambda_i^{(N)}$ and the sequences of quadratic (co-) variations of times $G^N(t),F^N(t),H^N(t)$ defined in Definition \ref{qvt} the following holds true:
\begin{enumerate}
\item[(i)] $G^N(t)\rightarrow G(t)~,~F^N(t)\rightarrow F(t)~,~H^N(t)\rightarrow H(t)$ as $N\rightarrow \infty$, where $G(t),F(t),H(t)$ are continuously differentiable functions on $[0,T]$.
\item[(ii)]For any null sequence $(h_N),\,h_N=\mathcal{O}\left(N^{-1}\right)$
\begin{subequations}\begin{align}\label{aqvtg}\frac{G^N(t+h_N)-G^N(t)}{h_N}\rightarrow G^{\prime}(t)\end{align} 
\begin{align}\label{aqvtf}\frac{F^N(t+h_N)-F^N(t)}{h_N}\rightarrow F^{\prime}(t)\end{align}
\begin{align}\label{aqvth}\frac{H^N(t+h_N)-H^N(t)}{h_N}\rightarrow H^{\prime}(t)\end{align}
\end{subequations}
uniformly on [0,T] as $N\rightarrow \infty$.
\end{enumerate}
\end{assump}
Assumption \ref{aqvt} is necessary to ensure that the sequence of variances of the estimator \eqref{HY} converges as $n,m\rightarrow\infty$. The derivative of the asymptotic quadratic variation of refresh times \eqref{aqvtg} will appear in the asymptotic variance of the discretization error $D_T^N$, since refresh times are (in general) not equidistant. For $\d T_i^{(N)}=T/N$ for all $i\in\{0,\ldots,N\}$, $G^{\prime}(t)=\1_{[0,T]}$ holds true.\\
The uniform convergence of the difference quotients defined by \eqref{aqvtf} and \eqref{aqvth} are necessary to ensure that the sequence of variances of $A_T^N$ converges as $N\rightarrow \infty$. The assumptions imposed by \eqref{aqvtg}-\eqref{aqvth} are weaker than assuming convergence of the joint sampling design of $\left(\mathcal{T}^{X,n},\mathcal{T}^{Y,m}\right)$ and are not very restrictive. They hold true whenever the sequences of sampling schemes tend to a certain state of asynchronicity or have a uniform  behaviour of non-synchronicity in the limit as $n,m\rightarrow\infty$. For homogeneous sampling schemes these (co-)variations of time converge to linear limiting functions.\\ 
The sequence of functions $F^N$ describe an interaction of interpolation steps between the two processes. In contrast, $H^N$ is defined to measure an impact of the in general non-zero correlations of next-tick and previous-tick interpolations to the same refresh time $T_i$, for each process separately.\\ \\
\textbf{Example:}\\
Consider the synchronous equidistant sampling schemes with $N=n=m$ and $t_i^{(n)}=\tau_{j}^{(n)}=i/n,i=0,\ldots,n$. The left-hand side of Figure \ref{figqvt} shows the quadratic (co)variations of time $G^N,F^N$ and $H^N$ for $N=30000$. $F^N$ and $H^N$ are identically zero since there are no asynchronous observations and because $T_i^{(N)}=i/N,t_i^{(n)}=\tau_i^{(n)},0\le i\le n$, interpolation steps are redundant and $A_T^N$ equals zero. The function $G^N$ is a step function that will tend to the identity on $[0,T]$ as $N\rightarrow\infty$.\\
Next, we consider a situation which originates from the complete synchronous equidistant one by shifting one time-scale half a time instant $1/2N$.
Then we have completely non-synchronous sampling schemes and we will call this situation intermeshed sampling.
In this case the synchronous approximation is still equidistant with instants $1/N$ and, hence, $G$ is the identity function. $F$ and $H$ are linear limiting functions with slope 1 and 1/4, respectively. Interpolations are carried out for all $1\le i\le N$ for the same process for which its first observation takes place after the first observation of the other process. All interpolation steps equal $1/2N$ and thus $H^{\prime}=1/4$ follows. Since for $H$ interpolated time instants $1/2N$ are multiplied with refresh time instants $1/N$ in \textbf{both} addends due to the specific structure, $F$ equals the identity on $[0,T]$. The functions $G^N,F^N,H^N$ for intermeshed sampling are illustrated in Figure \ref{figqvt} on the right-hand side.\\ \\
In the next section we will show that for an important special case, independent homogeneous Poisson sampling, \eqref{aqvtg}-\eqref{aqvth} are fulfilled when replacing deterministic convergence by convergence in probability. Furthermore, the stochastic limits $G^{\prime}(t),F^{\prime}(t),H^{\prime}(t)$ are calculated explicitly and are again constant on $[0,T]$. For data applications one can calculate easily empirical versions $\tilde G_{n,m}^{\prime}(t),\tilde F_{n,m}^{\prime}(t),\tilde H_{n,m}^{\prime}(t)$ of $G^{\prime}(t),F^{\prime}(t),H^{\prime}(t)$ and use those as estimators for \eqref{aqvtg}-\eqref{aqvth}.\\
The key result of this section is the following Theorem \ref{HYclt}. The detailed proof is postponed to the Appendix \ref{sec:7}. This result gives insight into the asymptotic distribution of the Hayashi-Yoshida estimator. It improves on the asymptotic normality result in \cite{hy2}, since the weak convergence is stable in the setting where we allow for random correlation, drift and volatility processes. The representation of the asymptotic variance using  \eqref{aqvtg}-\eqref{aqvth} differs from that in \cite{hy3}, where a similar stable convergence result is established, by the decomposition of the estimation error in \eqref{D} and \eqref{As} and the notion of (co-)variations of times. The latter provide helpful tools to describe the stylized facts and features of non-synchronous data and build the ground work for combined approaches for widespread generalizations and extensions of the underlying model. A very important one is the generalized multiscale estimator in \cite{bibinger2} when market microstructure noise effects are taken into account. 
\begin{theo}\label{HYclt}The estimation error of \eqref{HY} converges on the Assumptions \ref{eff}, \ref{grid} and \ref{aqvt} stably in law to a centred, mixed Gaussian distribution:
\begin{equation}\sqrt{N}\left(\sum_{i=1}^N\left(X_{g_i}-X_{l_i}\right)\left(Y_{\gamma_i}-Y_{\lambda_i}\right)-\left[ X\,,\,Y\right]_T\right)\stackrel{st}{\rightsquigarrow}\mathbf{N}\left(0\,,\,v_{D_T}+v_{A_T}\right)~,\end{equation} 
with the asymptotic variance
\begin{equation*}\hspace*{-.05cm}v_{D_T}\hspace*{-.05cm}+\hspace*{-.05cm}v_{A_T}\hspace*{-.05cm}=\hspace*{-.05cm}T\hspace*{-.1cm}\int_0^T\hspace*{-.1cm} G^{\prime}(t)\hspace*{-.05cm}\left(\sigma_t^X\sigma_t^Y\right)^2\hspace*{-.05cm}\left(\rho_t^2+1\right)dt+T\hspace*{-.1cm}\int_0^T\hspace*{-.1cm}\left(F^{\prime}(t)\hspace*{-.05cm}\left(\sigma_t^X\sigma_t^Y\right)^2\hspace*{-.05cm}+2
H^{\prime}(t)\hspace*{-.05cm}\left(\rho_t\sigma_t^X\sigma_t^Y\right)^2\hspace*{-.05cm}\right)dt\end{equation*}
where the two addends come from the asymptotic variances of $D_T^N$ and $A_T^N$, respectively.
\end{theo}

\section{Independent Poisson sampling\label{sec:5}}
In this section, we consider the model in which the sequences of observation times are supposed to be realizations of two homogeneous Poisson processes that are mutually independent and independent of the processes $ X$ and $Y$. \\
Thereto, let $\tilde n^{(n)}(t)$ and $\tilde m^{(n)}(t)$ be sequences of two independent homogeneous Poisson processes with parameters $Tn/\theta_1$ and $Tn/\theta_2$ ($n\in\N$), such that the waiting times between jumps of $\tilde n^{(n)}$ and $\tilde m^{(n)}$ are exponentially distributed with expectations $\E\left[\Delta t_i^{(n)}\right]=\theta_1/n$ and $\E\left[\Delta \tau_j^{(n)}\right]=\theta_2/n~,i\in\N,j\in\N$. Thus, $\tilde n^{(n)}(T)$ and $\tilde m^{(n)}(T)$ correspond to the sequences giving the numbers of observation times of $ X$ and $Y$ in the time span $[0,T]$. The increments of the sampling times of the closest synchronous approximation \eqref{partition} are maxima of the exponentially distributed waiting times and we obtain:
$$\Delta T_k^{(n)}\sim F(t)=1-\exp{\left(-\frac{tn}{\theta_1}\right)}-\exp{\left(-\frac{tn}{\theta_2}\right)}+\exp{\left(-tn\left(\frac{1}{\theta_1}+\frac{1}{\theta_2}\right)\right)}~,k\in\N\,.$$
Denote $\tilde N(T)^{(n)}=\max_{N\in\N}{\{\sum_{k=0}^{N}\Delta T_k^{(n)}\le T\}}$. We focus on the characteristics of the sampling schemes affecting the asymptotics of the synchronized realized covolatility estimator \eqref{HY}. In particular our interest is in the quadratic (co-)variations of times defined in Definition \ref{qvt}.
\begin{prop}In the independent homogeneous Poisson model for sampling schemes, it holds true that
\begin{subequations}
\begin{align}\label{pp1}G^N(t)\stackrel{p}{\longrightarrow} 2\left(1-\frac{2\theta_1^2\theta_2^2}{\theta_1^2\theta_2^2+(\theta_1^2+\theta_2^2)(\theta_1+\theta_2)^2}\right)\frac{t}{T}~~\left(=\frac{14}{9}\frac{t}{T}~~\text{if}~~\theta_1=\theta_2=\theta\right)~,\end{align}
\begin{align}\label{pp2}F^N(t)\stackrel{p}{\longrightarrow} \left(\frac{2\theta_1\theta_2}{(\theta_1^2+\theta_1\theta_2+\theta_2^2)}+\frac{4\theta_1^2\theta_2^2}{\left(\theta_1+\theta_2-\frac{\theta_1\theta_2}{\theta_1+\theta_2}\right)^2(\theta_1+\theta_2)^2}\right)\frac{t}{T}\\
\notag~~~~~~~~~~~~~~~~~~~~~~~~~~~~~~~~~~~~~~~~~~~~~~~~~~~~~~~~~~~~~~~~~~~~~~~~~~~~~~~~~~~~~~~~~\left(=\frac{10}{9}\frac{t}{T}~~\text{if}~~\theta_1=\theta_2=\theta\right)~,\end{align}
\begin{align}\label{pp3}H^N(t)\stackrel{p}{\longrightarrow} 2\left(\frac{1}{\left(\theta_1+\theta_2-\frac{\theta_1\theta_2}{\theta_1+\theta_2}\right)^2}\frac{\theta_1^2\theta_2^2}{(\theta_1+\theta_2)^2}\right)\frac{t}{T}~~\left(=\frac{2}{9}\frac{t}{T}~~\text{if}~~\theta_1=\theta_2=\theta\right)~.\end{align}
\end{subequations}
\end{prop}
\begin{proof}
Poisson processes are Markovian and the exponential distribution of the increments between arrival times is memoryless. Wald's identity ensures that $\E\left[\sum_{k=0}^{\tilde N(T)^{(n)}}\Delta T_k^{(n)}\right]=\E\left[\tilde N(T)^{(n)}\right]\E\left[\Delta T_1^{(n)}\right]$. For the proofs of these attributes and further information on properties of mutually independent homogeneous Poisson processes we refer interested readers to \cite{coxisham}.\\
First of all we ascertain that $t_i^{(n)}\ne \tau_j^{(n)}\,\forall\,(i,j)\in\{1,\ldots,\tilde n^{(n)}(T)\}\times\{1,\ldots,\tilde m^{(n)}(T)\}$ almost surely. For an arbitrarily fixed $i$, the expected values of next-tick, previous-tick and refresh time instants yield
\begin{align*}
\E\left[g_i^{(n)}-T_i^{(n)}\right]=\E\left[\left(g_i^{(n)}-T_i^{(n)}\right)\,\big|\,T_i^{(n)}=\gamma_i^{(n)}\right]\P\left(T_i^{(n)}=\gamma_i^{(n)}\right)=\frac{\theta_1}{n}\frac{\theta_2}{\theta_1+\theta_2}~,\end{align*}
\begin{align*}
\E\left[\gamma_i^{(n)}-T_i^{(n)}\right]=\frac{\theta_2}{n}\frac{\theta_1}{\theta_1+\theta_2}~,\end{align*}
\begin{align*}
\E\left[T_i^{(n)}-l_{i+1}^{(n)}\right]=\int_0^{\infty}y\frac{n}{\theta_2}e^{-\frac{yn}{\theta_2}}e^{-\frac{yn}{\theta_1}}\,dy=\frac{1}{n}\frac{\theta_1^2\theta_2}{(\theta_1+\theta_2)^2}~,\end{align*}
\begin{align*}
\E\left[T_i^{(n)}-\lambda_{i+1}^{(n)}\right]=\frac{1}{n}\frac{\theta_1\theta_2^2}{(\theta_1+\theta_2)^2}~,\end{align*}
\begin{align*}
\E\left[T_{i+1}^{(n)}-T_{i}^{(n)}\right]=\frac{\theta_1}{n}+\frac{\theta_2}{n}-\frac{1}{n}\frac{\theta_1\theta_2}{\theta_1+\theta_2}~.\end{align*}
The conditional expectations given that the $i$th refresh time $T_i^{(n)}=\gamma_i^{(n)}$ is an arrival time of $\tilde m^{(n)}$ yield 
$\E\left[T_{i+1}^{(n)}-T_i^{(n)}|T_i^{(n)}=\gamma_i^{(n)}\right]=\E\left[T_{i+1}^{(n)}-T_i^{(n)}\right]$ and $\E\left[T_{i}^{(n)}-l_{i+1}^{(n)}|T_i^{(n)}=\gamma_i^{(n)}\right]=\E\left[T_{i}^{(n)}-l_{i+1}^{(n)}\right]$, since the latter previous-tick interpolation is zero with probability 1 if $T_i^{(n)}\ne \gamma_i^{(n)}$. Only for $(T_i^{(n)}-\lambda_i^{(n)})$ the conditional expectation differs from the unconditional and can be calculated by further conditioning
\begin{align*}&\E\left[T_i^{(n)}-\lambda_i^{(n)}|T_i^{(n)}=\gamma_i^{(n)}\right]=\\
&~~\E\left[T_i^{(n)}-\lambda_i^{(n)}|T_i^{(n)}=\gamma_i^{(n)}\,,\,T_{i-1}^{(n)}=\lambda_i^{(n)}\right]\P\left(T_{i-1}^{(n)}=\lambda_i^{(n)}|T_{i}^{(n)}=\gamma_i^{(n)}\right)\\ &~~~~~+\E\left[T_i^{(n)}-\lambda_i^{(n)}|T_i^{(n)}=\gamma_i^{(n)}\,,\,T_{i-1}^{(n)}=l_i^{(n)}\right]\P\left(T_{i-1}^{(n)}=l_i^{(n)}|T_{i}^{(n)}=\gamma_i^{(n)}\right)\\
&=\left(\theta_1+\theta_2-\frac{\theta_1\theta_2}{\theta_1+\theta_2}\right)\frac{\theta_1}{\theta_1+\theta_2}+2\theta_1\frac{\theta_2}{\theta_1+\theta_2}~,\end{align*}
where the factor $2\theta_1$ in the second addend is simply the expectation of the waiting time for two jumps of $\tilde n$. Here, we have used some simplifying symmetry aspects, a rigorous proof using the density functions is obtained by calculation of
$$\E\left[T_{i}^{(n)}-\lambda_{i}^{(n)}\1_{\{\,T_i^{(n)}=\gamma_i^{(n)},T_{i-1}^{(n)}=\lambda_i^{(n)}\}}\right]=\hspace*{-0.05cm}\int_0^{\infty}\hspace*{-0.075cm}\int_x^{\infty}x\frac{n}{\theta_1}e^{-x\frac{n}{\theta_1}}e^{-y\frac{n}{\theta_2}}y\frac{n}{\theta_2}e^{-x\frac{n}{\theta_1}}e^{-y\frac{n}{\theta_2}}\,dx\,dy=\frac{2\theta_1\theta_2}{\theta_1+\theta_2}~.$$
The conditional expectations on $T_i^{(n)}=g_i^{(n)}$ are deduced analogously.
Since $\E\left[T_i^{(n)}-l_i^{(n)}\right]=$\\ $\E\left[T_i^{(n)}-T_{i-1}^{(n)}\right]+\E\left[T_{i-1}^{(n)}-l_i^{(n)}\right]$ and the (conditional) expectations of the products occurring in $G^N,F^N,$ $H^N$ equal the products of (conditional) expectations thanks to the memorylessness of exponential distributions, the latter results suffice to apply the law of large numbers to the empirical (co-)variations of times. For the asymptotics of $G^N(T),F^N(T)$ and $H^N(T)$, we conclude for the number of addends $\tilde N(T)^{(n)}$, that $\E\tilde N(T)^{(n)}=(T/\theta) n+\KLEINO(n)$ with $\theta=\theta_1+\theta_2-(\theta_1\theta_2)/(\theta_1+\theta_2)$ what follows from $\E\tilde N(T)^{(n)}\E\left[\Delta T_1^{(N)}\right]=T+\mathcal{O}_p(n^{-1})$ and $\var\left(\tilde N(T)^{(n)}\right)=\mathcal{O}(n^{-1})$ since
\begin{align*}\var\left(\sum_{k=0}^{\tilde N(T)^{(n)}}\Delta T_k^{(n)}\right)=\var\left(\tilde N(T)^{(n)}\right)\E\left[\left(\Delta T_1^{(n)}\right)^2\right]+\E\left[\tilde N(T)^{(n)}\right]\var\left(\Delta T_1^{(n)}\right)~.\end{align*}
The exact probability mass functions of the counting processes $\tilde N(t)^{(n)}$ associated with the maxima of the waiting times $\Delta t_i^{(n)},\Delta\tau_j^{(n)}$ have a quite complicated form, so that we only give the last two results on the expectation and the variance that are necessary for the proof of the proposition. \\
From the preceding conclusions, it follows that
\begin{align*}G^N(t)=\frac{\tilde N(T)^{(n)}}{T}\sum_{T_i^{(n)}\le t}\left(\Delta T_i^{(n)}\right)^2\stackrel{p}{\longrightarrow} \frac{n^2}{\theta^2}\left(\frac{2\theta_1^2}{n^2}+\frac{2\theta_2^2}{n^2}-2\left(\frac{\theta_1\theta_2}{(\theta_1+\theta_2)}\right)^2\frac{1}{n^2}\right)\frac{t}{T}~,\end{align*}
\begin{align*}F^N(t)&=\frac{\tilde N(T)^{(n)}}{T}\sum_{T_{i+1}^{(n)}\le t}(T_i^{(n)}-\lambda_i^{(n)})(g_i^{(n)}-T_i^{(n)})+\left(T_i^{(n)}-l_{i}^{(n)}\right)\left(\gamma_{i}^{(n)}-T_i^{(n)}\right)\\&~~~~~~~~~~~~~~~~~~~~~~~~~~~~~~+\d T_{i+1}^{(n)}\left(T_i^{(n)}-l_{i+1}^{(n)}\right)+\d T_{i+1}^{(n)}\left(T_i^{(n)}-\lambda_{i+1}^{(n)}\right)~\\
&~~~~~\stackrel{p}{\longrightarrow}\frac{t}{T\theta^2}\left(\frac{\theta_1\theta_2}{(\theta_1+\theta_2)}\left(2\theta_1+2\theta_2-2\frac{\theta_1\theta_2}{(\theta_1+\theta_2)}+\frac{2\theta_1\theta_2}{(\theta_1+\theta_2)}\right)\right. \\
&\left.~~~~~~~~~~~~~~~~~~~~~~~~~~ +\left(\theta_1+\theta_2-\frac{\theta_1\theta_2}{(\theta_1+\theta_2)}\right)\frac{\theta_1^2\theta_2+\theta_1\theta_2^2}{(\theta_1+\theta_2)^2}\right)~,\end{align*}
\begin{align*}
H^N(t)&=\frac{\tilde N(T)^{(n)}}{T}\sum_{T_{i+1}^{(n)}\le t}\left(T_i^{(n)}-l_{i+1}^{(n)}\right)\left(g_i^{(n)}-T_i^{(n)}\right)+\left(T_i^{(n)}-\lambda_{i+1}^{(n)}\right)\left(\gamma_i^{(n)}-T_i^{(n)}\right)\\
&\stackrel{p}{\longrightarrow}\frac{t}{T\theta^2}\frac{\theta_1^2\theta_2^2(\theta_1+\theta_2)}{(\theta_1+\theta_2)^3}~.\end{align*}
Inserting $\theta$ we obtain formulae \eqref{pp1}-\eqref{pp3}. In the evaluation of $G^N$ we have also used the second moment of $\Delta T_1^{(n)}$ which can be calculated using the above given distribution function.
\end{proof}
\begin{figure}[t]
\framebox{
\includegraphics[width=7.285cm]{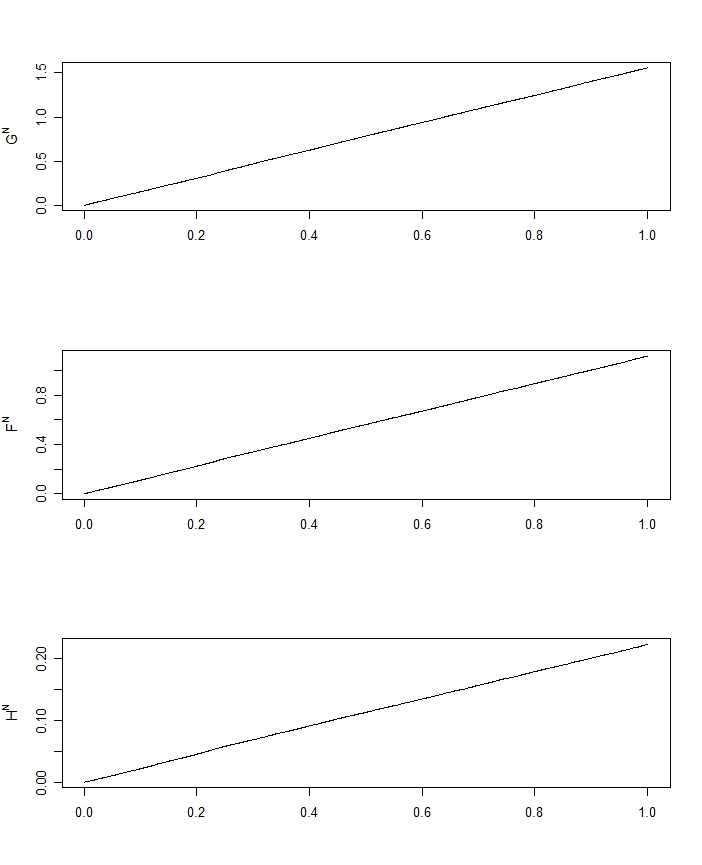}
\includegraphics[width=7.285cm]{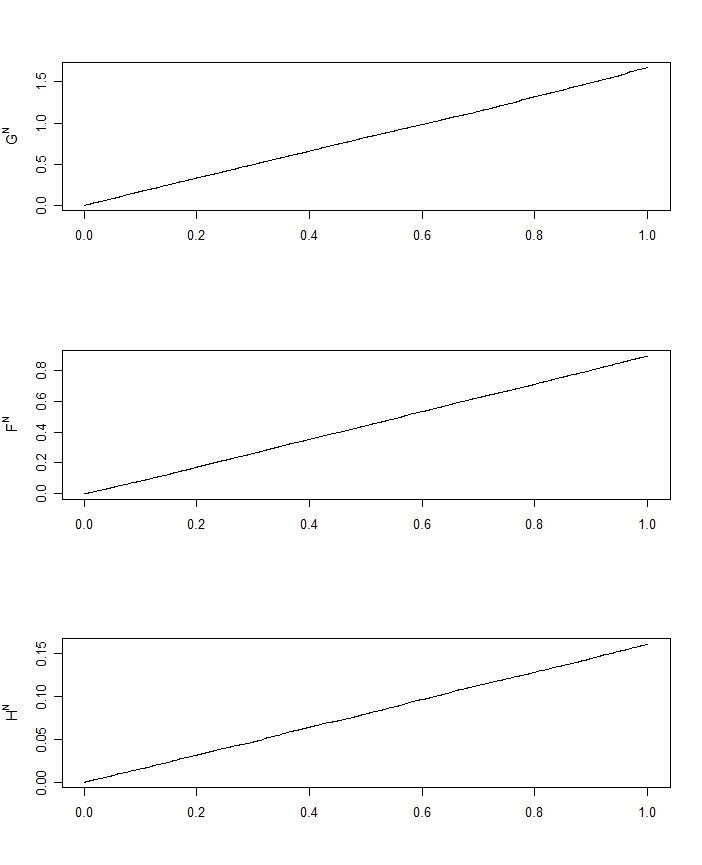}}
\caption{\label{figqvtpois}Quadratic (Co-)variations of times for homogeneous Poisson sampling.} 
\end{figure}
Figure \ref{figqvtpois} depitcs the quadratic (co-)variations of times for simulated mutually independent homogeneous Poisson processes. On the left-hand side both parameters have been set $\theta=1$ for $T=1$ and $n=30000$. The stochastic limits are linear increasing functions on $[0,1]$ with slope $14/9, 10/9, 2/9$ and $1/4$, respectively. On the right-hand side we see the (co-)variations of times for $T=1,n=30000,\theta_1=1,\theta_2=0.5$. Those tend in probability to linear limiting functions with slope $82/49,44/49,8/49$ and $2/9$, respectively.\\
In the model of non-synchronously observed It\^{o} processes $X$ and $Y$ which fulfill Assumption \ref{eff} and observation times following an independent Poisson sampling scheme of the above given form, we derive the following stable central limit theorem as special case of Theorem \ref{HYclt}:
\begin{cor}
\label{HYcltPP}The estimation error of the synchronized realized covolatility estimator \eqref{HY} converges on the Assumption \ref{eff} conditionally on the independent Poisson sampling scheme with $0<\theta_1<\infty$ and $0<\theta_2<\infty$ stably in law to a centred mixed Gaussian distribution:
\begin{equation}\sqrt{\tilde N(T)^{(n)}}\left(\sum_{i=0}^{\tilde N(T)^{(n)}}\left(X_{g_i^{(n)}}-X_{l_i^{(n)}}\right)\left(Y_{\gamma_i^{(n)}}-Y_{\lambda_i^{(n)}}\right)-\left[ X\,,\,Y\right]_T\right)\stackrel{st}{\rightsquigarrow}\mathbf{N}\left(0\,,\,v_T\right)~,\end{equation} 
with the asymptotic variance
\begin{equation*}v_T=2\int_0^T\left(\rho_t\sigma_t^X\sigma_t^Y\right)^2\,dt+\left(2\frac{\theta_1\theta_2}{\theta(\theta_1+\theta_2)}+1\right)\int_0^T\left(\sigma_t^X\sigma_t^Y\right)^2\end{equation*}
where the two addends come from the asymptotic variances of the discretization error $D_T^N$ of the closest synchronous approximation \eqref{D} and the additional error $A_T^N$ due to interpolations \eqref{As}, respectively, and $\theta=\theta_1+\theta_2-\frac{\theta_1\theta_2}{\theta_1+\theta_2}$.
\end{cor}
\begin{proof}It is a basic result in the theory of extreme values that for the supremum of $n$ i.\,i.\,d.\,exponentially distributed waiting times $\Delta T_i$ with $\E\Delta T_i=n^{-1}$, it holds true that $\sup_i{(\Delta T_i)}=\mathcal{O}_p\left(\log{(n)}/n\right)$. We refer to \cite{extrem} for a proof.
In the setting of mutually independent homogeneous Poisson processes with parameters $Tn/\theta_1$ and $Tn/\theta_2$, we conclude that $\sup_{i\in\{1,\ldots,\tilde N(T)^{(n)}\}}=\mathcal{O}_p\left(\log{\tilde N(T)^{(n)}}/\tilde N(T)^{(n)}\right)$. Hence, Assumption \ref{grid} holds for the sampling design where the orders of $\delta_n^X,\delta_m^Y$ hold in probability. Then all findings in the proofs of Propositions \ref{HYdisB} and \ref{HYasB} stay valid when we insert the (co-)variations of time deduced above in the limits of the variances.
\end{proof}
The stable convergence holds conditionally given the observation times, what means that endogenous observation times are not covered but Poisson sampling independent of the processes $X$ and $Y$.\\
The asymptotic variance of the mixed Gaussian limit is in line with the results by \cite{hy2} and \cite{hy3}. We remark that one has to pay attention to the proportionality to $\theta$ in the rate $\tilde N(T)^{(n)}$ when comparing the asymptotic variance to the one in \cite{hy3}. \\
From an applied point of view, the model considered in this section could be criticized for its flaw that sampling schemes of two correlated processes are modeled to follow two independent processes and for time homogeneity. Both seems to be rather unrealistic in financial time series. However, independent and homogeneous Poisson sampling times designs constitute the most commonly used model in this research area (cf.\,\cite{zhang}, \cite{hy} among others) because they are handy and allow for explicit calculations while the model is not too far away from the real world.

\section{Asymptotic variance estimation\label{sec:6}}
Finally, we state a consistent estimator for the asymptotic variance of the Hayashi-Yoshida estimator \eqref{HY} from Theorem \ref{HYclt}. Since in \cite{hy2} a central limit theorem for the case of deterministic correlation and volatility functions has been proved, the asymptotic variance is non-random in that setting. In a recent publication \cite{hy3}, in that the authors also generalize the asymptotic distribution result to a stable central limit theorem in the setting of random volatility and correlation functions, a consistent estimation method for the asymptotic variance is provided using kernel estimates. Our estimator differs from this method and we incorporate only one time transformed histogram-type estimator.
\begin{prop}\label{avarest}
Define the estimator
\begin{align*}\widehat{\AVAR}_{HY}\:=N\sum_{j=1}^{N-1}(X_{g_j}-X_{l_j})(Y_{\gamma_j}-Y_{\lambda_j})\left[(X_{g_j}-X_{l_j})(Y_{\gamma_j}-Y_{\lambda_j})\right. \\ \left. ~~~~~~+2(X_{g_{j+1}}-X_{l_{j+1}})(Y_{\gamma_{j+1}}-Y_{\lambda_{j+1}})\right]-3T\tilde I_1\end{align*}
with
\begin{align*}\tilde I_1\:=\sum_{j=1}^{K_N}\left(\frac{\widehat{\Delta\left[ X,Y\right]}_{G_j^N}^{HY}}{\Delta G_j^N}\right)^2\frac{G^N(T)}{K_N}\end{align*}
being a histogram-based estimator for $\int_0^T(\rho_t\sigma_t^X\sigma_t^Y)^2G^{\prime}(t)dt$. The estimators for the increase of the quadratic covariation on bins are Hayashi-Yoshida estimators of the type
\begin{align*}\widehat{\Delta\left[ X,Y\right]}_{G_j^N}^{HY}\:=\sum_{r\in[G_j^N,G_{j+1}^N)}(X_{g_r}-X_{l_r})(Y_{\gamma_r}-Y_{\lambda_r})~.\end{align*}
It holds true that
\begin{align*}\widehat{\AVAR}_{HY}\hspace*{-.1cm}\stackrel{p}{\longrightarrow}\hspace*{-.1cm}T\hspace*{-.1cm}\int_0^T\hspace*{-.1cm} G^{\prime}(t)\hspace*{-.05cm}\left(\sigma_t^X\sigma_t^Y\right)^2\hspace*{-.05cm}\left(\rho_t^2+1\right)dt+T\hspace*{-.1cm}\int_0^T\hspace*{-.1cm}\left(F^{\prime}(t)\hspace*{-.05cm}\left(\sigma_t^X\sigma_t^Y\right)^2\hspace*{-.05cm}dt+2
H^{\prime}(t)\hspace*{-.05cm}\left(\rho_t\sigma_t^X\sigma_t^Y\right)^2\hspace*{-.05cm}dt\right)\end{align*}
on the Assumptions \ref{eff}, \ref{grid} and \ref{aqvt}. Thus, we have on hand a consistent estimator for the asymptotic variance of the Hayashi-Yoshida estimator and the feasible stable central limit theorem
\begin{align}\label{hist}\frac{\widehat{\left[ X,Y\right]}_T^{(HY)}}{\sqrt{\widehat{\AVAR}_{HY}}}\stackrel{st}{\rightsquigarrow}\mathbf{N}(0,1)~.\end{align}
\end{prop}
For more motivation and details on the construction of histogram estimators for which bins are chosen equispaced according to a transformed timescale associated with a certain monotonic function, as the asymptotic quadratic variation of refresh times here, we refer to \cite{bibinger2}. Proposition \ref{avarest} is proved in Appendix \ref{sec:8}.
\renewcommand{\appendixname}{}
\appendix
 \section{Proof of Theorem \ref{HYclt}\label{sec:7}}
\subsection{Discretization error of the synchronous approximation}
\begin{prop}\label{HYdis}On the Assumptions \ref{eff}, \ref{grid} and \eqref{aqvtg} the discretization error of the closest synchronous approximation converges stably in law to a centred mixed Gaussian distribution:
\begin{align} \sqrt{\frac{N}{T}}D_T^{N} \stackrel{st}{\rightsquigarrow} \mathbf{N}\left(0\,,\,\int_0^T G^{\prime}(t)\left(\sigma_t^X\sigma_t^Y\right)^2(\rho_t^2+1)\,dt\right)~.
\end{align}
\end{prop}
\begin{proof}
In the proofs superscripts of the sampling times are frequently omitted to increase the readability.\\ 
First note that on Assumption \ref{eff}, by Girsanov's theorem we may without loss of generality further suppose that $\mu_t^X=\mu_t^Y=0$ identically since we have learned in Section \ref{sec:2} that stable convergence is commutative with measure change. Let $M_t$ and $L_t$ be the continuous martingales $L_t=\int_0^t\sigma_s^X\,dW_s^X~,~M_t=\int_0^t\sigma_s^Y\,dW_s^Y$ where $W^X,W^Y$ are two standard Brownian motions with quadratic covariation $\left[ W^X,W^Y\right]_t=\int_0^t\rho_s\sigma_s^X\sigma_s^Y\,ds$ and denote
$L_i=\int_0^{T_i}\sigma_s^X\,dW_s^X~,~M_i=\int_0^{T_i}\sigma_s^Y\,dW_s^Y$.
\begin{prop}\label{HYdisB}On the same Assumptions as in Proposition \ref{HYdis},
 the process $\mathcal{D}_t^N$ defined by
\begin{align*}\mathcal{D}_t^N\:=\sqrt{\frac{N}{T}}\sum_{T_{i}^{(N)}\le t}(L_{i}-L_{{i-1}})(M_{i}-M_{{i-1}})-\int_0^t\rho_s\sigma_s^X\sigma_s^Yds\end{align*}
for $0\le t\le T$ converges as $N\rightarrow\infty$ stably in law:
\begin{align}\mathcal{D}_t^N\stackrel{st}{\rightsquigarrow}\int_0^t\sqrt{v_{\mathcal{D}_s}}dW^{\bot}_s\end{align}
where $W^{\bot}$ is a Brownian motion independent of $\mathcal{F}$ and
\begin{align}v_{\mathcal{D}_s}=G^{\prime}(s)(\sigma_s^X\sigma_s^Y)^2(\rho_s^2+1)~.\end{align}
\end{prop}
\begin{proof}
We will prove this stable convergence of the process associated with the transformed discretization error by application of Jacod's stable limit Theorem \ref{jacod}. It is also possible to use the discrete-time version of this Theorem from Corollary \ref{djac} which we apply in the next subsection. \\
Using the definition of the quadratic covariation process of martingales or integration by parts formula, we find an illustration of the discretization error by a sum of stochastic integrals and an asymptotically negligible term:
\allowdisplaybreaks[3]{
\begin{align*}&\sum_{T_{i}^{(N)}\le t}\left(L_{T_i}-L_{T_{i-1}}\right)\left(M_{T_i}-M_{T_{i-1}}\right)=\sum_{T_{i}^{(N)}\le t}\left(L_i-L_{i-1}\right)\left(M_i-M_{i-1}\right)\\
&=\sum_{T_{i}^{(N)}\le t}\left(L_i M_i-L_i M_{i-1}-M_i L_{i-1}+L_{i-1}M_{i-1}\right)\\
&=\sum_{T_{i}^{(N)}\le t}\Big(\int_{T_{i-1}}^{T_i}L_sdM_s+\int_{T_{i-1}}^{T_i}M_sdL_s+\d\left[ L,M\right]_{T_i}\\ &\phantom{~=~~~~~~~~~}-M_{i-1}(L_i-L_{i-1})-L_{i-1}(M_i-M_{i-1})\Big)\\
&=\left[ L,M\right]_{\mathfrak{T}(t)}-\left[ L,M\right]_{T_0}+\sum_{T_{i}^{(N)}\le t}\left(\int_{T_{i-1}}^{T_i}(L_s-L_{i-1})dM_s+\int_{T_{i-1}}^{T_i}(M_s-M_{i-1})dL_s\right)
\end{align*}}
where we denote $\mathfrak{T}(t)\:=\max_i{(T_i^{(N)}\le t)}$.\\ 
Thus, we obtain 
$$\sqrt{\frac{N}{T}}\mathcal{D}_t^N=\sqrt{\frac{N}{T}}\sum_{T_{i}^{(N)}\le t}\left(\int_{T_{i-1}}^{T_i}(L_s-L_{i-1})dM_s+\int_{T_{i-1}}^{T_i}(M_s-M_{i-1})dL_s\right)+\KLEINO_p\left(\sqrt{\frac{N}{T}}\right)~,$$
since $\left[ L,M\right]_{\mathfrak{T}(t)}-\left[ L,M\right]_{T_0}=\left[ L,M\right]_t+\KLEINO_p(1)$.
Consider the centred continuous martingale
\begin{align*}\phi_{\tau}^{(N)}\:=\sqrt{\frac{N}{T}}\Big(\sum_{T_{i}^{(N)}\le t}\left(\int_{T_{i-1}}^{T_i}(L_s-L_{i-1})dM_s+\int_{T_{i-1}}^{T_i}(M_s-M_{i-1})dL_s\right)\hspace*{2.6cm}\\\hspace*{2.6cm}+\int_{\mathfrak{T}(t)}^{\tau}(L_s-L_{\mathfrak{T}(t)})dM_s+\int_{\mathfrak{T}(t)}^{\tau}(M_s-M_{\mathfrak{T}(t)})dL_s\Big)~,~\tau \in [\mathfrak{T}(t),t]~.\end{align*}
We calculate the corresponding quadratic variation process at time $t$:
{\allowdisplaybreaks[2]{
\begin{align*}\left[ \phi^{(N)}\right]_t&=\frac{N}{T}\left[\sum_{T_{i}^{(N)}\le t}\left(\int_{T_{i-1}}^{T_i}(L_s-L_{i-1})dM_s+\int_{T_{i-1}}^{T_i}(M_s-M_{i-1})dL_s\right)\right]_t\hspace*{-.15cm}+\left[\phi^{(N)}\right]_t-\hspace*{-.05cm}\left[\phi^{(N)}\right]_{\mathfrak{T}(t)}\\
&=\frac{N}{T}\sum_{T_{i}^{(N)}\le t}\Big(\int_{T_{i-1}}^{T_i}(L_s-L_{i-1})^2 d\left[ M\right]_s+\int_{T_{i-1}}^{T_i}(M_s-M_{i-1})^2d\left[ L\right]_s\\
&~~~~~+2\int_{T_{i-1}}^{T_i}(L_s-L_{i-1})(M_s-M_{i-1})d\left[ M,L\right]_s\Big)+\left[\phi^{(N)}\right]_t-\left[\phi^{(N)}\right]_{\mathfrak{T}(t)}\\
&{\underset{\text{\tiny{(Lemma \ref{HYdisBlem})}}}{=}}\frac{N}{T}\sum_{T_{i}^{(N)}\le t}\Big(\int_{T_{i-1}}^{T_i}\left[ L-L_{i-1}\right]_s d\left[ M\right]_s+\int_{T_{i-1}}^{T_i}\left[ M-M_{i-1}\right]_s d\left[ L\right]_s\\
&~~~~~+2\int_{T_{i-1}}^{T_i}\left[ L-L_{i-1}\right]_s\left[ M-M_{i-1}\right]_s d\left[ M,L\right]_s\Big)+\KLEINO_p(1)\\
&=\frac{N}{T}\sum_{i=1}^N\Big(\int_{T_{i-1}}^{T_i}\,d\left(\left[ L-L_{i-1}\right]_s \left[ M-M_{i-1}\right]_s\right)\\ &~~~~~+2\int_{T_{i-1}}^{T_i}\left[ L-L_{i-1}\right]_s\left[ M-M_{i-1}\right]_s d\left[ M,L\right]_s\Big)+\KLEINO_p(1)\\
&=\frac{N}{T}\sum_{T_{i}^{(N)}\le t}\left(\int_{T_{i-1}}^{T_i}\left(\sigma_s^X\right)^2\,ds\int_{T_{i-1}}^{T_i}\left(\sigma_s^Y\right)^2\,ds\right)+\left(\int_{T_{i-1}}^{T_i}\rho_s\sigma_s^X\sigma_s^Y\,ds\right)^2+\KLEINO_p(1)\\
&=\frac{N}{T}\sum_{T_{i}^{(N)}\le t}\left(\left((\overline{\rho\sigma^X\sigma^Y})_i\right)^2\left(\d T_i\right)^2+\left((\overline{\sigma^X})_i(\overline{\sigma^Y})_i\right)^2\left(\d T_i\right)^2\right)+\KLEINO_p(1)\\
&=\sum_{T_{i}^{(N)}\le t}\left(\frac{\left(G^{(N)}(T_{i})-G^{(N)}(T_{i-1})\right)}{\d T_i}\left(\left(\sigma^X_{T_{i-1}}\sigma^Y_{T_{i-1}}\right)^2\left(1+{\rho_{T_{i-1}}}^2\right)\right)\,\d T_i\right)+\KLEINO_p(1)\\
&\stackrel{p}{\longrightarrow} \int_0^tG^{\prime}(s)(\rho_s^2+1)(\sigma_s^X\sigma_s^Y)^2\,ds~.
\end{align*}}}
In this calculation we have used integration by parts and the change of variables Theorem for the integrals with quadratic covariation integrators that are of finite variation. The second last equality is an application of the mean value theorem (the volatility and the correlation processes are continuous and thus also bounded on compact sets) where the constants $(\overline{\sigma^X})_i$, $(\overline{\sigma^Y})_i$ and $(\overline{\rho\sigma^X\sigma^Y})_i$ come from. The Riemann sum converges and with Definition \ref{qvt} and Assumption \ref{aqvt} this yields the convergence in probability of the quadratic variation to $\int_0^tG^{\prime}(s)(\rho_s^2+1)(\sigma_s^X\sigma_s^Y)^2\,ds=\int_0^t v_{\mathcal{D}_s}$.
The third equality above is proved in: 
\begin{lem}\label{HYdisBlem}It holds true that the approximation error terms
\begin{subequations}
\begin{align}\label{subdis1}\sum_{T_{i}^{(N)}\le t}\hspace*{-0.125cm}\left(\int_{T_{i-1}}^{T_i}(L_s-L_{i-1})^2d\left[ M-M_{i-1}\right]_s-\int_{T_{i-1}}^{T_i}\left[ L-L_{i-1}\right]_sd\left[ M-M_{i-1}\right]_s\right)\end{align}
\begin{align}\label{subdis2}\sum_{T_{i}^{(N)}\le t}\hspace*{-0.05cm}\left(\int_{T_{i-1}}^{T_i}(M_s-M_{i-1})^2d\left[ L-L_{i-1}\right]_s-\int_{T_{i-1}}^{T_i}\left[ M-M_{i-1}\right]_s d\left[ L-L_{i-1}\right]_s\right)\end{align}
\begin{align}\label{subdis3}\sum_{T_{i}^{(N)}\le t}\hspace*{-0.05cm}\left(\int_{T_{i-1}}^{T_i}\hspace*{-0.1cm}(M_s-M_{i-1})(L_s-L_{i-1})d\left[ M,L\right]_s-\hspace*{-0.075cm}\int_{T_{i-1}}^{T_i}\hspace*{-0.1cm}\left[ M-M_{i-1},L-L_{i-1}\right]_s d\left[ M,L\right]_s\hspace*{-0.05cm}\right)\end{align}
\begin{align}\label{subdis4}\frac{N}{T}\sum_{T_{i}^{(n)}\le t}(\Delta T_i)^2\left((\overline{\rho\sigma^X\sigma^Y})_i^2+(\overline{\sigma^X})_i^2(\overline{\sigma^Y})_i^2-\left(\left(\rho_{T_{i-1}}\sigma^X_{T_{i-1}}\sigma^Y_{T_{i-1}}\right)^2+\left(\sigma^X_{T_{i-1}}\sigma^Y_{T_{i-1}}\right)^2\right)\right)\end{align}
\end{subequations}
converge to zero in probability.
\end{lem}
\begin{proof}
The proofs for \eqref{subdis1} and \eqref{subdis2} are completely analogous and we restrict ourselves to prove it for \eqref{subdis1}. By It\^{o}'s formula
$$(L_s-L_{i-1})^2=2\int_{T_{i-1}}^s(L_r-L_{i-1})dL_r+\left[ L-L_{i-1}\right]_s$$
holds. The left-hand side of \eqref{subdis1} equals
\begin{align*}&\sum_{T_{i}^{(N)}\le t}\left(\int_{T_{i-1}}^{T_i}\left(2\int_{T_{i-1}}^s(L_r-L_{i-1})dL_r\right)d\left[ M-M_{i-1}\right]_r\right)\\
&=\sum_{T_{i}^{(N)}\le t}\left(2\int_{T_{i-1}}^{T_i}(L_s-L_{i-1})(\left[ M-M_{i-1}\right]_{T_i})\,dL_s-2\int_{T_{i-1}}^{T_i}(L_s-L_{i-1})(\left[ M-M_{i-1}\right]_{s})\,dL_s\right)\end{align*}
by application of the integration by parts formula in the way
$$Z_{T_i}\left[ M-M_{i-1}\right]_{T_i}=\int_0^{T_i}Z_td\left[ M-M_{i-1}\right]_t+\int_0^{T_i}\left[ M-M_{i-1}\right]_t dZ_t$$
with $Z_t\:=\int_{T_{i-1}}^t2(L_s-L_{i-1})dL_s$ for $T_{i-1}\le t\le T$ to the addends. Therefore, we can write the left-hand side of \eqref{subdis1} in the way $\mathcal{M}_1^{(N)}+\mathcal{M}_2^{(N)}$ with two centred continuous martingales $\mathcal{M}_1^{(N)}, \mathcal{M}_2^{(N)}$ defined in the fashion of $\phi^{(N)}$ above and calculate the quadratic covariation processes at time $t$:
\begin{align*}\left[ \mathcal{M}_2^{(N)}\right]_t&=4\sum_{T_{i}^{(N)}\le t}\left(\int_{T_{i-1}}^{T_i}(L_s-L_{i-1})^2(\left[ M-M_{i-1}\right]_s)^2d\left[ L\right]_s\right)+\KLEINO_p(1)\\
&\le 4\max_{i}{\sup_{s\in(T_{i-1},T_i]}{\hspace*{-0.075cm}(L_s-L_{i-1})^2}}\hspace*{-0.06cm}\max_{i}{\sup_{s\in(T_{i-1},T_i]}{\hspace*{-0.075cm}\left[ M-M_{i-1}\right]_s^2}}\hspace*{-.1cm}\sum_{T_{i}^{(N)}\le t}\int_{T_{i-1}}^{T_i}d\left[ L\right]_s+\KLEINO_p(1)~.\end{align*}
The first addend is up to a logarithmic factor $\mathcal{O}_p(\delta_N^{3})$ and hence $\mathcal{M}_2^{(N)}=\KLEINO_p(1)$ on Assumption \ref{grid}. That $\mathcal{M}_1^{(N)}=\KLEINO_p(1)$ is proved analogously. This implies that \eqref{subdis1} is $\KLEINO_p(1)$. \\
The strategy of the proof for \eqref{subdis3} follows the same approach, starting with the equation
\begin{align*}(L_s-L_{i-1})(M_s-M_{i-1})&=\int_{T_{i-1}}^{s}(L_r-L_{i-1})d(M-M_{i-1})_r+\int_{T_{i-1}}^s(M_r-M_{i-1})d(L-L_{i-1})_r\\ &~~+\left[ L-L_{i-1},M-M_{i-1}\right]_s\end{align*}
and applying integration by parts as above with $Z_t=\int_{T_{i-1}}^t(L_s-L_{i-1})d(M-M_{i-1})_s+\int_{T_{i-1}}^t(M_s-M_{i-1})d(L-L_{i-1})_s$ for $T_{i-1}\le t\le T_i$. \\
We complete the proof of the convergence of the quadratic variation with the proof for \eqref{subdis4}. Denote $\overline{(\rho\sigma^X\sigma^Y)}_i^2=\widetilde {(\sigma^X)^2}_i\cdot \widetilde {(\sigma^Y)^2}_i\cdot \widetilde {(\rho)^2}_i$ to distinguish between the values from the application of the mean value theorems to the two different addends. An upper bound of the left-hand side of \eqref{subdis4} can be found by elementary algebra and the triangle inequality for the absolute value:
\begin{align*}&
\frac{N}{T}\sum_{T_{i}^{(N)}\le t}(\Delta T_i)^2\left((\widetilde{\rho}\widetilde{\sigma^X}\widetilde{\sigma^Y})_i^2+(\overline{\sigma^X})_i^2(\overline{\sigma^Y})_i^2-\left(\left(\rho_{T_{i-1}}\sigma^X_{T_{i-1}}\sigma^Y_{T_{i-1}}\right)^2+\left(\sigma^X_{T_{i-1}}\sigma^Y_{T_{i-1}}\right)^2\right)\right)\\
&\le \frac{N}{T}\sum_{T_{i}^{(N)}\le t}(\Delta T_i)^2\left(\widetilde {(\sigma^X)^2}_i\widetilde {(\sigma^Y)^2}_i\left|\widetilde {(\rho)^2}_i-\rho_{T_{i-1}}^2\right|+\widetilde {(\sigma^Y)^2}_i\rho_{T_{i-1}}^2\left|\widetilde {(\sigma^X)^2}_i-(\sigma^X_{T_{i-1}})^2\right|\right. \\
&~~~~\left. +\rho_{T_{i-1}}^2(\sigma_{T_{i-1}}^X)^2\left|\widetilde {(\sigma^Y)^2}_i-(\sigma^Y_{T_{i-1}})^2\right|+(\overline{\sigma^Y})_i^2\left|(\overline{\sigma^X})_i^2-(\sigma^X_{T_{i-1}})^2\right|+(\sigma^X_{T_{i-1}})^2\left|(\overline{\sigma^Y})_i^2-(\sigma^Y_{T_{i-1}})^2\right|\right)\\
&=\KLEINO_p(1)~.
\end{align*}
\end{proof}
The martingales $\phi^{(N)}$ can be written for every $N$ as time-changed Brownian motions $B^{(DDS,N)}_{\left[ \phi^{(N)}\right]_t}=\phi^{(N)}_t$ by the Dambis-Dubins-Schwarz theorem. The sequence of martingales $\phi^{(N)}$ or associated time-changed Dambis-Dubins-Schwarz Brownian motions converges weakly to a limiting Brownian motion by the asymptotic Knight-theorem. The limiting Brownian motion will be defined on an orthogonal extension of the original probability space. To obtain the stable convergence result, we apply Jacod's Theorem \ref{jacod} and thus, we are left to verify conditions \eqref{conditionsstable1} and \eqref{conditionsstable2}.\\
Consider the quadratic covariation process of $\phi^{(N)}$ and the reference martingale $L$
$$\left[ L,\phi^{(N)}\right]_t=\sqrt{\frac{N}{T}}\sum_{T_{i}^{(N)}\le t}\left(\int_{T_{i-1}}^{T_i}(L_s-L_{i-1})d\left[ M,L\right]_s+\int_{T_{i-1}}^{T_i}(M_s-M_{i-1})d\left[ L\right]_s\right)+\KLEINO_p(1)~.$$
The term of smaller order than 1 in probability comes from the increment of the covariation process on $[\mathfrak{T}(t),t]$. As before, this equality holds true for all $t$, since for $t<T_1$ the covariation is $\KLEINO_p(1)$.
Integration by parts yields:
\begin{align*}\left[ L,\phi^{(N)}\right]_t&=\sqrt{\frac{N}{T}}\sum_{T_{i}^{(N)}\le t}\left[(\left[ M,L\right]_{T_i}-\left[ M,L\right]_{T_{i-1}})(L_i-L_{i-1})-\int_{T_{i-1}}^{T_i}\left[ M,L\right]_sd(L_s-L_{i-1})\right. \\ &\left.\hspace*{3cm} +(\left[ L\right]_{T_i}-\left[ L\right]_{T_{i-1}})(M_i-M_{i-1})-\int_{T_{i-1}}^{T_i}\left[ L\right]_sd(M_s-M_{i-1})\right]~.\end{align*}
It remains to show that this term converges to zero in probability. The term is centred and using It\^{o} isometry we find the following upper bound for the second moment:
\begin{align*}\E\left[\left(\left[ L,\phi^{(N)}\right]_t\right)^2\right]&\le 2 \frac{N}{T}\E\left[\sum_{T_i^{(N)}\le t}\left((\left[ M,L\right]_{T_i}-\left[ M,L\right]_{T_{i-1}})^2(L_i-L_{i-1})^2\right.\right. \\ &\hspace*{1.5cm}\left.\left. +(\left[ L\right]_{T_i}-\left[ L\right]_{T_{i-1}})^2(M_i-M_{i-1})^2\right)\right. \\ &\left.~~+\max_{i\in\{1,\ldots,N\}}{\sup_{s\in(T_{i-1},T_i]}{(\left[ M,L\right]_{s}-\left[ M,L\right]_{T_{i-1}})^2}}\sum_i\int_{T_{i-1}}^{T_i}d\left[ L-L_{i-1}\right]_t \right. \\ &\hspace*{1.5cm}\left. +\max_{i\in\{1,\ldots,N\}}{\sup_{s\in(T_{i-1},T_i]}{(\left[ L\right]_{s}-\left[ L\right]_{T_{i-1}})^2}}\sum_i\int_{T_{i-1}}^{T_i}d\left[ M-M_{i-1}\right]_t\right] \\
&=\mathcal{O}\left(N\delta_N^2\right)~~.\end{align*}
The term is bounded by a constant times $N\delta_N^2$ since squared increments, cross products of increments and increments of the quadratic (co-)variations of $L$ and $M$ over time instants $\Delta T_i^{(N)}$ are bounded by $\Delta T_i^{(N)}$ times a constant. To sums with products of time instants we can apply Hölder's inequality with the supremum norm to obtain upper bounds. There are at most order $\delta_N^{-1}$ time instants $\Delta T_i^{(N)}$ of order $\sup_i{\Delta T_i^{(N)}}=\delta_N$ since $\sum_i\Delta T_i^{(N)}\le T$ and the time span $T$ is fixed.\\
Hence, $\left[ L,\phi^{(N)}\right]_t=\KLEINO_p(1)~~~\forall t \,\in[0,T]$. With the same strategy $\left[ M,\phi^{(N)}\right]_t=\KLEINO_p(1)~~~\forall t\, \in[0,T]$ can be shown. \\
For every bounded $\mathcal{F}_t$-martingale $L^{\bot}$ satisfying $\left[ L,L^{\bot}\right]\equiv 0$ the covariation
\begin{align*}\left[ L^{\bot},\phi^{(N)}\right]_t=\frac{N}{T}\sum_{T_{i}^{(N)}\le t}\left(\int_{T_{i-1}}^{T_i}(L_s^{\bot}-L_{i-1}^{\bot})d\left[ M,L^{\bot}\right]_s\right)+\KLEINO_p(1)=\KLEINO_p(1)\end{align*}
converges to zero. The same holds true for every bounded $\mathcal{F}_t$-martingale orthogonal to $M$. 
Applying Theorem \ref{jacod}, we deduce that Proposition \ref{HYdisB} holds true.
\end{proof}
Proposition \ref{HYdis} is a direct consequence of the stronger result in Proposition \ref{HYdisB} since for $t=T$ the marginal distribution is a mixed normal distribution which is independent of $\mathcal{F}$. The stable convergence assures that the convergence also holds under the original probability measure and non-zero drift terms with the same asymptotic law.
\end{proof}

\subsection{Error due to non-synchronicity} 
\begin{prop}\label{HYas}Let Assumptions \ref{eff}, \ref{grid} and \eqref{aqvtf}-\eqref{aqvth} from Assumption \ref{aqvt} be satisfied. The error $A_T^N$ due to the lack of synchronicity converges stably in law to a centred mixed Gaussian distribution:
\begin{align}\sqrt{\frac{N}{T}}A_T^N\stackrel{st}{\rightsquigarrow}\mathbf{N}\left(0,v_{A_T}\right)~,\end{align}
with asymptotic variance
\begin{align}v_{A_T}=\int_0^TF^{\prime}(t)\left(\sigma_t^X\sigma_t^Y\right)^2dt+\int_0^T 2H^{\prime}(t)\left(\rho_t\sigma_t^X\sigma_t^Y\right)^2dt~.\end{align}
\end{prop}
\begin{proof}
First, we write the $i$th increments occurring as factors in the addends of the estimator \eqref{HY} as the sum of the next-tick interpolation at $T_i$, the increments $\Delta X_{T_i}=X_{T_i}-X_{T_{i-1}}$ and $\Delta Y_{T_i}=Y_{T_i}-Y_{T_{i-1}}$, respectively, and the previous-tick interpolation at $T_{i-1}$ and multiply out the addends.
\begin{align*}\widehat{\left[ X,Y\right]}_T&=\sum_{i=1}^N\left(X_{g_i}\hspace*{-0.05cm}-\hspace*{-0.05cm}X_{T_i}\hspace*{-0.05cm}+\hspace*{-0.05cm}X_{T_i}\hspace*{-0.05cm}-\hspace*{-0.05cm}X_{T_{i-1}}\hspace*{-0.05cm}+\hspace*{-0.05cm}X_{T_{i-1}}\hspace*{-0.05cm}-\hspace*{-0.05cm}X_{l_i}\right)\left(Y_{\gamma_i}\hspace*{-0.05cm}-Y_{T_i}\hspace*{-0.05cm}+\hspace*{-0.05cm}Y_{T_i}\hspace*{-0.05cm}-\hspace*{-0.05cm}Y_{T_{i-1}}\hspace*{-0.05cm}+\hspace*{-0.05cm}Y_{T_{i-1}}\hspace*{-0.05cm}-\hspace*{-0.05cm}Y_{\lambda_i}\right)\\
&=\sum_{i=1}^N\left((X_{g_i}\hspace*{-0.05cm}-\hspace*{-0.05cm}X_{T_i})\Delta Y_{T_i}+(Y_{\gamma_i}\hspace*{-0.05cm}-\hspace*{-0.05cm}Y_{T_i})\Delta X_{T_i}+(X_{T_{i-1}}\hspace*{-0.05cm}-\hspace*{-0.05cm}X_{l_i})\Delta Y_{T_i}+(Y_{T_{i-1}}-Y_{\lambda_i})\Delta X_{T_i}\right.\\
&~~\left.+(X_{g_i}-X_{T_{i}})(Y_{T_{i-1}}-Y_{\lambda_i})+(Y_{\gamma_i}-Y_{T_i})(X_{T_{i-1}}-X_{l_i})\right)+D_T^N
\end{align*}
The indicator functions in \eqref{As} have been dropped since the corresponding addends are zero if the indicator functions were zero. Since at least one of the next-tick interpolation errors is zero and as well one of the previous-tick interpolation errors, too, two addends, namely the products of next-tick interpolation errors and the product of previous-tick interpolation errors, equal zero. Thus, the error due to asynchronicity can be written as the sum of the remaining six terms (where at least another three equal zero in each addend).
We conclude, that the error $A_T^N$ can be expressed in the following way:
\begin{align*}
A_T^N&=\sum_{i=1}^{N-1}\left((X_{g_i}-X_{T_i})(Y_{T_i}-Y_{\lambda_i})+(Y_{\gamma_{i}}-Y_{T_i})(X_{T_i}-X_{l_{i}})\right.\\
&~~~~~~~~~~~~~\left.(X_{T_{i+1}}-X_{T_i})(Y_{T_i}-Y_{\lambda_{i+1}})+(X_{T_i}-X_{l_{i+1}})(Y_{T_{i+1}}-Y_{T_i})\right)+\KLEINO_p(1)~.\end{align*}
In this equality an index shift has been applied to the partial sum of previous-tick interpolated errors multiplied with $\Delta X_{T_i}$ and $\Delta Y_{T_i}$, respectively, leading to the structure that in the $i$th addend the factors contain next- and previous-tick interpolated errors to the same $T_i$. The $\KLEINO_p(1)$-term emerges from end-effects when shifting the original sum.\\
In the last illustration of $A_T^N$ consecutive addends of the sum are uncorrelated in contrast to the non-shifted illustration. The reason is that, if without loss of generality $\gamma_i=T_i$ holds, $(X_{g_i}-X_{T_i})\Delta Y_{T_i}$ and $(X_{T_i}-X_{l_{i+1}})\Delta Y_{T_{i+1}}$ have in general a non-zero correlation whereas $(X_{g_i}-X_{T_i})\Delta Y_{T_i}$ and $(X_{T_{i-1}}-X_{l_{i}})\Delta Y_{T_{i}}$ are uncorrelated. Furthermore, the fact that $\gamma_i=T_i\Rightarrow \lambda_{i+1}=T_i$ assures that the addends in the last illustration of $A_T^N$ are uncorrelated. Roughly speaking we capture correlation between subsequent addends of the outer sum and transfer it into additional correlation in the inner sum.\\ 
As in the foregoing proof of Proposition \ref{HYdis}, it is sufficient to prove the stable convergence result for the zero-drift case. We denote, as before, the corresponding transformed processes $L_t=\int_0^t\sigma_s^X dW_s^X$ and $M_t=\int_0^t\sigma_s^Y dW_s^Y$.\\
Consider the sum
\begin{align}\notag\mathcal{A}_t^N=\sum_{T_{i+1}^{(N)}\le t}\Delta A_i^{N}\:&=\sqrt{\frac{N}{T}}\sum_{T_{i+1}^{(N)}\le t}\left((L_{g_i}-L_{T_i})(M_{T_i}-M_{\lambda_i})+(M_{\gamma_i}-M_{T_i})(L_{T_i}-L_{l_i})
\right.\\
&~~\left.+(L_{T_i}-L_{l_{i+1}})(M_{T_{i+1}}-M_{T_i})+(M_{T_i}-M_{\lambda_{i+1}})(L_{T_{i+1}}-L_{T_i})\right)\end{align}
for fixed $0\le t\le T$.
\begin{prop}\label{HYasB} Assume the same conditions as in Proposition \ref{HYas}. For fixed $0\le t\le T$ the transformed error due to non-synchronicity $\mathcal{A}_t^N$ is the endpoint of a discrete, centred, square-integrable martingale with respect to the filtration
$\mathcal{F}_{i,N}\:=\mathcal{F}_{T_{i+1}^{(N)}}$. The process $\mathcal{A}_t^N$ converges as $N\rightarrow \infty$ stably in law:
\begin{align}\mathcal{A}_t^N\stackrel{st}{\rightsquigarrow}\mathcal{A}_t=\int_0^t \sqrt{v_{\mathcal{A}_s}}dW^{\bot}_s\end{align}
where $W^{\bot}$ is a Brownian motion independent of $\mathcal{F}$ and 
\begin{align}v_{\mathcal{A}_s}= F^{\prime}(s)\left(\sigma_s^X\sigma_s^Y\right)^2+ 2H^{\prime}(s)\left(\rho_s\sigma_s^X\sigma_s^Y\right)^2~.\end{align}
\end{prop}
\begin{proof}
The expectation of the absolute value of the sum is bounded for all $t\in[0,T]$ and $\Delta A_i^{N},\,i=0,\ldots,N$ are $\mathcal{F}_{i,N}=\mathcal{F}_{T_{i+1}^{(N)}}$-measurable.
Since 
\begin{align*}\E\left[\Delta A_i^N|\mathcal{F}_{i-1,N}\right]&=\E\left[\Delta A_{i}^N|\mathcal{F}_{T_i^{(N)}}\right]\\
&=\E\left[(L_{g_i}-L_{T_i})(M_{T_i}-M_{\lambda_i})+(M_{\gamma_i}-M_{T_i})(L_{T_i}-L_{l_i})\right. \\
& ~~~~~~\left.+(L_{T_i}-L_{l_{i+1}})\Delta M_{T_{i+1}}+(M_{T_i}-M_{\lambda_{i+1}})\Delta L_{T_{i+1}}|\mathcal{F}_{T_i^{(N)}}\right]
\\ &=\E\left[L_{g_i}-L_{T_i}\right](M_{T_i}-M_{\lambda_i})+\E\left[M_{\gamma_i}-M_{T_i}\right](L_{T_i}-L_{l_i})\\ &~~~~~~+(L_{T_i}-L_{l_{i+1}})\E\left[\Delta M_{T_{i+1}}\right]+(M_{T_i}-M_{\lambda_{i+1}})\E\left[\Delta L_{T_{i+1}}\right]=0\end{align*}
for the conditional expectation of the increments holds, $\mathcal{A}_t^N$ is the endpoint of a $\mathcal{F}_{i,N}$-martingale.\\
The stable weak convergence to a limiting Brownian motion is proven with Corollary \ref{djac} to Jacod's Theorem \ref{jacod}.\\
First, we verify the conditional Lindeberg condition that is implied by the stronger conditional Lyapunov condition. It is sufficient to proof the following:
\begin{lem}\label{HYasLy}The sum of the conditional fourth moments of the martingale increments $A_i^N$ converges to zero in probability:
\begin{align*}\E\left[\sum_{T_{i+1}^{(N)}\le t}\left(\Delta A_i^N\right)^4\Big|\mathcal{F}_{i-1,N}\right]=\KLEINO_p(1)~.\end{align*}\end{lem}
\begin{proof}
Throughout the proof $C$ denotes a generic constant that does not depend on $N$. We consider different addends of the fourth conditional moments consecutively. The sum of conditional fourth moments incorporates addends of the following types:
\begin{itemize}
\item fourth-order moments: $$\frac{N^2}{T^2}\sum_{T_{i+1}^{(N)}\le t}\E\left[(L_{g_i}-L_{T_i})^4\right](M_{T_i}-M_{\lambda_i})^4~,$$
\item second-order moments: $$\frac{N^2}{T^2}\sum_{T_{i+1}^{(N)}\le t}\E\left[(L_{g_i}-L_{T_i})^2(\Delta M_{T_{i+1}})^2\right](L_{T_i}-L_{l_{i+1}})^2(M_{T_i}-M_{\lambda_i})^2~,$$
\item third- and first-order moments: $$\frac{N^2}{T^2}\sum_{T_{i+1}^{(N)}\le t}4(M_{T_i}-M_{\lambda_i})^3(L_{T_i}-L_{l_{i+1}})^3\E\left[\Delta M_{T_{i+1}}(L_{g_i}-L_{T_i})\right]~.$$
\end{itemize}
For the partial sum with addends of the first type an application of the Burkholder-Davis-Gundy (BDG) inequalities yields
\begin{align*}&\frac{N^2}{T^2}\sum_{T_{i+1}^{(N)}\le t}\E\left[(L_{g_i}-L_{T_i})^4\right](M_{T_i}-M_{\lambda_i})^4\\
\le &C\,\frac{N^2}{T^2}\sum_{T_{i+1}^{(N)}\le t}\E\left[\left(\int_{T_i}^{g_i}(\sigma_s^X)^2ds\right)^2\right](M_{T_i}-M_{\lambda_i})^4\\
\le &C\,\frac{N^2}{T^2}\sup_{s\in[0,T]}(\sigma_s^X)^2\sum_{T_{i+1}^{(N)}\le t}(M_{T_i}-M_{\lambda_i})^4(g_i-T_i)^2\le \mathcal{O}_p\left(N\delta_N^2\right)=\KLEINO_p(1)~.\end{align*}
The last inequality can be deduced by the result that the convergence $(N/(3T))\sum_i(\Delta M_{T_i})^4\rightarrow \int_0^t(\sigma_s^Y)^4ds$ holds almost surely as $N\rightarrow\infty$ for the so-called realized quarticity (\cite{bn3}) and that $(g_i-T_i)\le \delta_N$. Without the result about the convergence of the realized quarticity, the asymptotic order in probability can be derived by the convergence to zero of the expectation of the above sum and calculating the second moment that is bounded from above by a constant times $N^{4}\delta_N^7$.\\
For the partial sum including addends that incorporate second-order moments we obtain an upper bound by application of the Cauchy-Schwarz inequality and the BDG inequalities:
{\allowdisplaybreaks[2]{
\begin{align*}&\frac{N^2}{T^2}\sum_{T_{i+1}^{(N)}\le t}6\,\E\left[(L_{g_i}-L_{T_i})^2(\Delta M_{T_{i+1}})^2\right](L_{T_i}-L_{l_{i+1}})^2(M_{T_i}-M_{\lambda_i})^2\\
\le &\frac{N^2}{T^2}\sum_{T_{i+1}^{(N)}\le t}6\,\sqrt{\E\left[(L_{g_i}-L_{T_i})^4\right]}\sqrt{\E\left[(\Delta M_{T_{i+1}})^4\right]}(L_{T_i}-L_{l_{i+1}})^2(M_{T_i}-M_{\lambda_i})^2\\
\le &C\,\frac{N^2}{T^2}\sum_{T_{i+1}^{(N)}\le t}6\left(\E\left[\left(\int_{T_i}^{g_i}(\sigma_s^X)^2ds\right)^2\right]\E\left[\left(\int_{T_i}^{T_{i+1}}(\sigma_s^Y)^2ds\right)^2\right]\right)^{\frac{1}{2}}(L_{T_i}-L_{l_{i+1}})^2(M_{T_i}-M_{\lambda_i})^2\\ &=\KLEINO_p(1)~.\end{align*}
The stochastic order follows, since the term has the expectation
\begin{align*}&C\,\frac{N^2}{T^2}\hspace*{-0.1cm}\sum_{T_{i+1}^{(N)}\le t}\hspace*{-0.1cm}6\hspace*{-0.05cm}\left(\hspace*{-0.05cm}\E\left[\left(\int_{T_i}^{g_i}(\sigma_s^X)^2ds\right)^2\right]\hspace*{-0.05cm}\E\left[\left(\int_{T_i}^{T_{i+1}}(\sigma_s^Y)^2ds\right)^2\right]\right)^{\frac{1}{2}}\hspace*{-0.1cm}\E\left[(L_{T_i}-L_{l_{i+1}})^2(M_{T_i}-M_{\lambda_i})^2\right]\\
&\le C\,\frac{N^2}{T^2}\hspace*{-0.1cm}\sum_{T_{i+1}^{(N)}\le t}\hspace*{-0.1cm}6\hspace*{-0.05cm}\left(\hspace*{-0.05cm}\E\left(\int_{T_i}^{g_i}(\sigma_s^X)^2ds\hspace*{-0.05cm}\right)^2\hspace*{-0.1cm}\E\left(\int_{T_i}^{T_{i+1}}(\sigma_s^Y)^2ds\hspace*{-0.05cm}\right)^2\hspace*{-0.1cm}\E\left(\int_{l_{i+1}}^{T_i}(\sigma_s^X)^2ds\hspace*{-0.05cm}\right)^2\hspace*{-0.1cm}\E\left(\int_{\lambda_i}^{T_{i}}(\sigma_s^Y)^2ds\hspace*{-0.05cm}\right)^2\right)^{\frac{1}{2}}\\ &\le C\,N^2\delta_N^3=\KLEINO(1)~,\end{align*}}}
where again the Cauchy-Schwarz and BDG inequalities have been applied. The variance is bounded from above by a constant times $N^{4}\delta_N^7$, what can be shown by a similar calculation where thanks to the fact that $T_i=\gamma_i\Rightarrow \lambda_{i+1}=T_i$ the addends are uncorrelated and the variance of the sum equals the sum of variances.\\
We treat the third type of addends occurring in the sum of conditional fourth moments in the same way. It\^{o} isometry yields
\begin{align*} \frac{N^2}{T^2}\sum_{T_{i+1}^{(N)}\le t}4(M_{T_i}-M_{\lambda_i})^3(L_{T_i}-L_{l_{i+1}})^3\E\left[\Delta M_{T_{i+1}}(L_{g_i}-L_{T_i})\right]\\
=\frac{N^2}{T^2}\sum_{T_{i+1}^{(N)}\le t}4(M_{T_i}-M_{\lambda_i})^3(L_{T_i}-L_{l_{i+1}})^3\E\left[\int_{T_i}^{g_i}\rho_s\sigma_s^X\sigma_s^Yds\right]~.\end{align*}
This term has expectation
\begin{align*}&\frac{N^2}{T^2}\sum_{T_{i+1}^{(N)}\le t}4\E\left[(M_{T_i}-M_{\lambda_i})^3(L_{T_i}-L_{l_{i+1}})^3\right]\E\left[\int_{T_i}^{g_i}\rho_s\sigma_s^X\sigma_s^Yds\right]\\
& \le \frac{N^2}{T^2}\sum_{T_{i+1}^{(N)}\le t}4\sqrt{\E\left[(M_{T_i}-M_{\lambda_i})^6\right]\E\left[(L_{T_i}-L_{l_{i+1}})^6\right]}\E\left[\int_{T_i}^{g_i}\rho_s\sigma_s^X\sigma_s^Yds\right]\displaybreak[1]\\
& \le C\,\frac{N^2}{T^2}\sum_{T_{i+1}^{(N)}\le t}4\left(\E\left[\left(\int_{\lambda_i}^{T_i}(\sigma_s^Y)^2ds\right)^3\right]\E\left[\left(\int_{l_{i+1}}^{T_i}(\sigma_s^X)^2ds\right)^3\right]\right)^{\nicefrac{1}{2}}\E\left[\int_{T_i}^{g_i}\rho_s\sigma_s^X\sigma_s^Yds\right]\\
&\le C N^2\delta_N^3=\KLEINO(1)~,\end{align*}
and an analogous calculation as before yields that the variance is of order $N^4\delta_N^7$.\\
Thereby, the sum converges to zero in probability.
\end{proof}
Next, we consider the sum of conditional variances of the increments of the discrete martingale.
\begin{lem}\label{HYasV}\begin{align}\E\left[\sum_{T_{i+1}^{(N)}\le t}\left(\Delta A_i^N\right)^2\Big|\mathcal{F}_{T_{i}^{(N)}}\right]\stackrel{p}{\rightarrow}\int_0^tF^{\prime}(s)\left(\sigma_s^X\sigma_s^Y\right)^2ds+\int_0^t2H^{\prime}(s)\left(\rho_s\sigma_s^X\sigma_s^Y\right)^2ds~.\end{align}
\end{lem}
It holds true that
\begin{proof}
{\allowdisplaybreaks[2]{
\begin{align*}&\E\left[\sum_{T_{i+1}^{(N)}\le t}\left(\Delta A_i^N\right)^2\Big|\mathcal{F}_{T_{i}^{(N)}}\right]\\
&=\frac{N}{T}\hspace*{-0.05cm}\sum_{T_{i+1}^{(N)}\le t}\hspace*{-0.05cm}\E\left[(L_{g_i}\hspace*{-0.05cm}-\hspace*{-0.05cm}L_{T_i})^2(M_{T_i}\hspace*{-0.05cm}-\hspace*{-0.05cm}M_{\lambda_i})^2\hspace*{-0.05cm}+\hspace*{-0.05cm}(M_{\gamma_i}\hspace*{-0.05cm}-\hspace*{-0.05cm}M_{T_i})^2(L_{T_i}\hspace*{-0.05cm}-\hspace*{-0.05cm}L_{l_i})^2+(L_{T_i}\hspace*{-0.05cm}-\hspace*{-0.05cm}L_{l_{i+1}})^2(\Delta M_{T_{i+1}})^2 \right. \\ &\left. ~~~~~~~~~~~~+(M_{T_i}-M_{\lambda_{i+1}})^2(\Delta L_{T_{i+1}})^2+2(L_{g_i}-L_{T_i})(M_{T_i}-M_{\lambda_i})(L_{T_i}-L_{l_{i+1}})\Delta M_{T_{i+1}}\right.\\ &\left.~~~~~~~~~~~~+2(M_{\gamma_i}-M_{T_i})(L_{T_i}-L_{l_i})(M_{T_i}-M_{l_{i+1}})\Delta L_{T_{i+1}}\Big|\mathcal{F}_{T_i^{(N)}}\right]\displaybreak[0] \\
&{\underset{\text{\tiny{(It\^{o} isometry)}}}{=}}\frac{N}{T}\sum_{T_{i+1}^{(N)}\le t}\left(\E\left[\int_{T_i}^{g_i}(\sigma_s^X)^2ds\right](M_{T_i}-M_{\lambda_i})^2+\E\left[\int_{T_i}^{\gamma_i}(\sigma_s^Y)^2ds\right](L_{T_i}-L_{l_i})^2 \right.
\\ &\left.~~~~~~~~~~~~+(L_{T_i}-L_{l_{i+1}})^2\E\left[\int_{T_i}^{T_{i+1}}(\sigma_s^Y)^2ds\right]+(M_{T_i}-M_{\lambda_{i+1}})^2\E\left[\int_{T_i}^{T_{i+1}}(\sigma_s^X)^2ds\right]\right.
\\ &\left. ~~~~~~~~~~~~+2(M_{T_i}-M_{\lambda_i})(L_{T_i}-L_{l_{i+1}})\E\left[\int_{T_i}^{g_i}\rho_s\sigma_s^X\sigma_s^Yds\right]\right.
\\ &\left.~~~~~~~~~~~~+2(L_{T_i}-L_{l_i})(M_{T_i}-M_{l_{i+1}})\E\left[\int_{T_i}^{\gamma_i}\rho_s\sigma_s^X\sigma_s^Yds\right]\right)\displaybreak[0] \\
&{\underset{\text{\tiny{(Lemma \ref{h1})}}}{=}}\frac{N}{T}\sum_{T_{i+1}^{(N)}\le t}\left(\E\left[\int_{T_i}^{g_i}(\sigma_s^X)^2ds\right]\int_{\lambda_i}^{T_i}(\sigma_s^Y)^2ds+\E\left[\int_{T_i}^{\gamma_i}(\sigma_s^Y)^2ds\right]\int_{l_i}^{T_i}(\sigma_s^X)^2ds \right.
\\ &\left.~~~~~~~~~~~~+\int_{l_{i+1}}^{T_i}(\sigma_s^X)^2ds\E\left[\int_{T_i}^{T_{i+1}}(\sigma_s^Y)^2ds\right]+\int_{\lambda_{i+1}}^{T_i}(\sigma_s^Y)^2ds\E\left[\int_{T_i}^{T_{i+1}}(\sigma_s^X)^2ds\right]\right. 
\\ &\left. ~~~~~~~~~~~~+\int_{l_{i+1}}^{T_i}\rho_s\sigma_s^X\sigma_s^Yds\,\E\left[\int_{T_i}^{g_i}\rho_s\sigma_s^X\sigma_s^Yds\right]+\int_{l_{i+1}}^{T_i}\rho_s\sigma_s^X\sigma_s^Yds\,\E\left[\int_{T_i}^{\gamma_i}\rho_s\sigma_s^X\sigma_s^Yds\right]\right)+\KLEINO_p(1)\displaybreak[0]\\
&{\underset{\text{\tiny{(Lemma \ref{h1})}}}{=}}\frac{N}{T}\sum_{T_{i+1}^{(N)}\le t}\left(\int_{T_i}^{g_i}(\sigma_s^X)^2ds\int_{\lambda_i}^{T_i}(\sigma_s^Y)^2ds+\int_{T_i}^{\gamma_i}(\sigma_s^Y)^2ds\int_{l_i}^{T_i}(\sigma_s^X)^2ds \right. 
\\ &\left.~~~~~~~~~~~~+\int_{l_{i+1}}^{T_i}(\sigma_s^X)^2ds\int_{T_i}^{T_{i+1}}(\sigma_s^Y)^2ds+\int_{\lambda_{i+1}}^{T_i}(\sigma_s^Y)^2ds\int_{T_i}^{T_{i+1}}(\sigma_s^X)^2ds\right.
\\ &\left. ~~~~~~~~~~~~+2\int_{l_{i+1}}^{T_i}\rho_s\sigma_s^X\sigma_s^Yds\int_{T_i}^{g_i}\rho_s\sigma_s^X\sigma_s^Yds+2\int_{l_{i+1}}^{T_i}\rho_s\sigma_s^X\sigma_s^Yds\int_{T_i}^{\gamma_i}\rho_s\sigma_s^X\sigma_s^Yds\right)+\KLEINO_p(1) \displaybreak[0]\\
&{\underset{\text{\tiny{(Lemma \ref{h2})}}}{=}}\frac{N}{T}\sum_{T_{i+1}^{(N)}\le t}\left((\sigma_{T_i}^X\sigma_{T_i}^Y)^2\left((T_i-\lambda_i)(g_i-T_i)+(\gamma_i-T_i)(T_i-l_i)+(T_i-l_{i+1})\Delta T_{i+1}\right.\right.\\
&\left.\left.~~~~~~~~~~~~\hspace*{-0.25cm}+(T_i-\lambda_{i+1})\Delta T_{i+1}\right)\hspace*{-0.05cm}+\hspace*{-0.05cm} (\rho_{T_i}\sigma_{T_i}^X\sigma_{T_i}^Y)^2\hspace*{-0.05cm}\left(2(T_i-l_{i+1})(g_i-T_i)\hspace*{-0.05cm}+\hspace*{-0.05cm}2(T_i-\lambda_{i+1})(\gamma_i-T_i)\right)\right)\hspace*{-0.05cm}+\hspace*{-0.05cm}\KLEINO_p(1)\displaybreak[0]\\
&=\sum_{T_{i+1}^{(N)}\le t}\frac{F(T_{i+1})-F(T_i)}{T_{i+1}-T_i}(\sigma_{T_i}^X\sigma_{T_i}^Y)^2\Delta T_{i+1}+2\frac{H(T_{i+1})-H(T_i)}{T_{i+1}-T_i}(\rho_{T_i}\sigma_{T_i}^X\sigma_{T_i}^Y)^2\Delta T_{i+1}+\KLEINO_p(1)\\
&\stackrel{p}{\longrightarrow}\int_0^tF^{\prime}(s)\left(\sigma_s^X\sigma_s^Y\right)^2ds+\int_0^t2H^{\prime}(s)\left(\rho_s\sigma_s^X\sigma_s^Y\right)^2ds~.
\end{align*}}}
In the last step we have involved Definition \ref{qvt}. The Riemann sum converges on the Assumption \ref{aqvt} (in particular \eqref{aqvtf} and \eqref{aqvth}) in probability as $N\rightarrow \infty$ to the expression $\int_0^t v_{\mathcal{A}_s}ds$ with $v_{\mathcal{A}_s}$ given in Proposition \ref{HYasB}.\\
The detailed proofs of the approximations are postponed in the following two lemmas.
\begin{lem}\label{h1}On the assumptions as before, the following equations hold true:
{\allowdisplaybreaks[2]{
\begin{align*}\frac{N}{T}\sum_{T_{i+1}^{(N)}\le t}\left((M_{T_i}-M_{\lambda_i})^2-\int_{\lambda_i}^{T_i}(\sigma_s^Y)^2ds\right)\E\left[\int_{T_i}^{g_i}(\sigma_s^X)^2ds\right]=\KLEINO_p(1)~,\end{align*}
\begin{align*}\frac{N}{T}\sum_{T_{i+1}^{(N)}\le t}\int_{\lambda_i}^{T_i}(\sigma_s^Y)^2ds\left(\E\left[\int_{T_i}^{g_i}(\sigma_s^X)^2ds\right]-\int_{T_i}^{g_i}(\sigma_s^X)^2ds\right)=\KLEINO_p(1)~,\end{align*}
\begin{align*}\frac{N}{T}\sum_{T_{i+1}^{(N)}\le t}\left((L_{T_i}-L_{l_i})^2-\int_{l_i}^{T_i}(\sigma_s^X)^2ds\right)\E\left[\int_{T_i}^{\gamma_i}(\sigma_s^Y)^2ds\right]=\KLEINO_p(1)~,\end{align*}
\begin{align*}\frac{N}{T}\sum_{T_{i+1}^{(N)}\le t}\int_{l_i}^{T_i}(\sigma_s^X)^2ds\left(\E\left[\int_{T_i}^{\gamma_i}(\sigma_s^Y)^2ds\right]-\int_{T_i}^{\gamma_i}(\sigma_s^Y)^2ds\right)=\KLEINO_p(1)~,\end{align*}
\begin{align*}\frac{N}{T}\sum_{T_{i+1}^{(N)}\le t}\left((M_{T_i}-M_{\lambda_{i+1}})^2-\int_{\lambda_{i+1}}^{T_i}(\sigma_s^Y)^2ds\right)\E\left[\int_{T_i}^{T_{i+1}}(\sigma_s^X)^2ds\right]=\KLEINO_p(1)~,\end{align*}
\begin{align*}\frac{N}{T}\sum_{T_{i+1}^{(N)}\le t}\int_{\lambda_{i+1}}^{T_i}(\sigma_s^Y)^2ds\left(\E\left[\int_{T_i}^{T_{i+1}}(\sigma_s^X)^2ds\right]-\int_{T_i}^{T_{i+1}}(\sigma_s^X)^2ds\right)=\KLEINO_p(1)~,\end{align*}
\begin{align*}\frac{N}{T}\sum_{T_{i+1}^{(N)}\le t}\left((L_{T_i}-L_{l_{i+1}})^2-\int_{l_{i+1}}^{T_i}(\sigma_s^X)^2ds\right)\E\left[\int_{T_i}^{T_{i+1}}(\sigma_s^Y)^2ds\right]=\KLEINO_p(1)~,\end{align*}
\begin{align*}\frac{N}{T}\sum_{T_{i+1}^{(N)}\le t}\int_{l_{i+1}}^{T_i}(\sigma_s^X)^2ds\left(\E\left[\int_{T_i}^{T_{i+1}}(\sigma_s^Y)^2ds\right]-\int_{T_i}^{T_{i+1}}(\sigma_s^Y)^2ds\right)=\KLEINO_p(1)~,\end{align*}
\begin{align*}\frac{N}{T}\sum_{T_{i+1}^{(N)}\le t}\left((M_{T_i}-M_{\lambda_{i}})(L_{T_i}-L_{l_{i+1}})-\int_{l_{i+1}}^{T_i}\rho_s\sigma_s^X\sigma_s^Yds\right)\E\left[\int_{T_i}^{g_i}\rho_s\sigma_s^X\sigma_s^Yds\right]=\KLEINO_p(1)~,\end{align*}
\begin{align*}\frac{N}{T}\sum_{T_{i+1}^{(N)}\le t}\int_{l_{i+1}}^{T_i}\rho_s\sigma_s^X\sigma_s^Yds\left(\E\left[\int_{T_i}^{g_i}\rho_s\sigma_s^X\sigma_s^Yds\right]-\int_{T_i}^{g_i}\rho_s\sigma_s^X\sigma_s^Yds\right)=\KLEINO_p(1)~,\end{align*}
\begin{align*}\frac{N}{T}\sum_{T_{i+1}^{(N)}\le t}\left((L_{T_i}-L_{l_{i}})(M_{T_i}-L_{\lambda_{i+1}})-\int_{\lambda_{i+1}}^{T_i}\rho_s\sigma_s^X\sigma_s^Yds\right)\E\left[\int_{T_i}^{\gamma_i}\rho_s\sigma_s^X\sigma_s^Yds\right]=\KLEINO_p(1)~,\end{align*}
\begin{align*}\frac{N}{T}\sum_{T_{i+1}^{(N)}\le t}\int_{\lambda_{i+1}}^{T_i}\rho_s\sigma_s^X\sigma_s^Yds\left(\E\left[\int_{T_i}^{\gamma_i}\rho_s\sigma_s^X\sigma_s^Yds\right]-\int_{T_i}^{\gamma_i}\rho_s\sigma_s^X\sigma_s^Yds\right)=\KLEINO_p(1)~.\end{align*}}}
\end{lem}
\begin{proof}
We restrict ourselves to the proof of the first two equalities, since all other terms can shown to converge to zero in probability in an analogous way. The left-hand side of the first equality has an expectation equal to zero which can be concluded directly by It\^{o} isometry:
$$\E\left[\frac{N}{T}\sum_{T_{i+1}^{(N)}\le t}\left((M_{T_i}-M_{\lambda_i})^2-\int_{\lambda_i}^{T_i}(\sigma_s^Y)^2ds\right)\E\left[\int_{T_i}^{g_i}(\sigma_s^X)^2ds\right]\right]=0~.$$
In order to derive the stochastic order of the term, consider the second moment:
\begin{align*}&\E\left[\left(\frac{N}{T}\sum_{T_{i+1}^{(N)}\le t}\left((M_{T_i}-M_{\lambda_i})^2-\int_{\lambda_i}^{T_i}(\sigma_s^Y)^2ds\right)\E\left[\int_{T_i}^{g_i}(\sigma_s^X)^2ds\right]\right)^2\right]\\
&=\frac{N^2}{T^2}\sum_{T_{i+1}^{(N)}\le t}\E\left[(M_{T_i}-M_{\lambda_i})^4-2(M_{T_i}-M_{\lambda_i})^2\int_{\lambda_i}^{T_i}(\sigma_s^Y)^2ds+\left(\int_{\lambda_i}^{T_i}(\sigma_s^Y)^2ds\right)^2\right]\\ &~~~~~~\times\left(\E\left[\int_{T_i}^{g_i}(\sigma_s^X)^2ds\right]\right)^2
=\KLEINO(1)~,
\end{align*}
where the asymptotic order is deduced by It\^{o} isometry and the BDG inequalities. Since the error induced by this term in the approximation of the conditional variance before is centred and has a variance converging to zero as $N\rightarrow\infty$, the error is asymptotically negligible.\\
In the second equality we consider the error when the expected increment of the quadratic variation of $X$ over the next-tick interpolated time interval is substituted by the integral itself. We proceed as before for the first approximation. Since
$$\E\left[\frac{N}{T}\sum_{T_{i+1}^{(N)}\le t}\int_{\lambda_i}^{T_i}(\sigma_s^Y)^2ds\left(\E\left[\int_{T_i}^{g_i}(\sigma_s^X)^2ds\right]-\int_{T_i}^{g_i}(\sigma_s^X)^2ds\right)\right]=0$$
and
\begin{align*}
&\E\left[\left(\frac{N}{T}\sum_{T_{i+1}^{(N)}\le t}\int_{\lambda_i}^{T_i}(\sigma_s^Y)^2ds\left(\E\left[\int_{T_i}^{g_i}(\sigma_s^X)^2ds\right]-\int_{T_i}^{g_i}(\sigma_s^X)^2ds\right)\right)^2\right]\\
&=\frac{N^2}{T^2}\sum_{T_{i+1}^{(N)}\le t}\var\left(\int_{\lambda_i}^{T_i}(\sigma_s^Y)^2ds\right)\E\left[\left(\int_{\lambda_i}^{T_i}(\sigma_s^Y)^2ds\right)^2\right]=\KLEINO(1)~,\end{align*}
the approximation error is asymptotically negligible. The fact that $\gamma_i=T_i\Rightarrow \lambda_{i+1}=T_i$ has been used that guarantees that the addends of the sum are uncorrelated.
\end{proof}
Lemma \ref{h1} has been applied in the second and third equality in the evaluation of the sum of conditional variances and the proof of Lemma \ref{HYasV} is completed by the following
\begin{lem}\label{h2}On the same assumptions as before, the following equation holds true 
\begin{align*}\frac{N}{T}\sum_{T_{i+1}^{(N)}\le t}\left(\int_{T_i}^{g_i}(\sigma_s^X)^2ds\int_{\lambda_i}^{T_i}(\sigma_s^Y)^2ds-(\sigma_{T_i}^X\sigma_{T_i}^Y)^2(T_i-\lambda_i)(g_i-T_i)\right)=\KLEINO_p(1)\end{align*}
and analogously the errors in the five other addends converge to zero in probability when replacing the product of increments of quadratic (co-)variations by the values of $\rho_{T_i},\sigma^X_{T_i},\sigma^Y_{T_i}$ multiplied with the corresponding times increments.
\end{lem}
\begin{proof}
We prove the equality explicitly given in the lemma. The five remaining terms can be handled by the same strategy. By an application of the mean value theorem, elementary algebra and the triangle inequality for the absolute values, we deduce that
{\allowdisplaybreaks[2]{
\begin{align*}&\frac{N}{T}\sum_{T_{i+1}^{(N)}\le t}\left(\int_{T_i}^{g_i}(\sigma_s^X)^2ds\int_{\lambda_i}^{T_i}(\sigma_s^Y)^2ds-(\sigma_{T_i}^X\sigma_{T_i}^Y)^2(T_i-\lambda_i)(g_i-T_i)\right)\\
&=\frac{N}{T}\sum_{T_{i+1}^{(N)}\le t}\left((\sigma_{\varsigma_i}^X\sigma_{\xi_i}^Y)^2-(\sigma_{T_i}^X\sigma_{T_i}^Y)^2\right)(T_i-\lambda_i)(g_i-T_i)\\
&\le \frac{N}{T}\sum_{T_{i+1}^{(N)}\le t}\left|(\sigma_{\varsigma_i}^X\sigma_{\xi_i}^Y)^2-(\sigma_{T_i}^X\sigma_{T_i}^Y)^2\right|(T_i-\lambda_i)(g_i-T_i)\\
&\le\frac{N}{T}\sum_{T_{i+1}^{(N)}\le t}\left(\left((\sigma_{\xi_i}^Y)^2\left|(\sigma_{\varsigma_i}^X)^2-(\sigma_{T_i}^X)^2\right|+(\sigma_{T_i}^X)^2\left|(\sigma_{\xi_i}^Y)^2-(\sigma_{T_i}^Y)^2\right|\right)(T_i-\lambda_i)(g_i-T_i)\right)\\
&\le C\,\frac{N}{T}\sum_{T_{i+1}^{(N)}\le t}\left(\left(\sup_{s\in[\lambda_i,T_i]}{\left|(\sigma_{s}^X)^2-(\sigma_{T_i}^X)^2\right|}+\sup_{s\in[T_i,g_i]}{\left|(\sigma_{s}^Y)^2-(\sigma_{T_i}^Y)^2\right|}\right)(T_i-\lambda_i)(g_i-T_i)\right)\\
&=\KLEINO_p(1)
\end{align*}}}
holds on Assumption \ref{eff}.
\end{proof}
\end{proof}
To prove the stability of the convergence in Proposition \ref{HYasB}, we show in the following that the discrete covariations of $\mathcal{A}_t^N$ with the $\mathcal{F}$-generating underlying martingales $L$ and $M$ converge to zero in probability. 
\begin{lem}\label{HYasC1}
\begin{align*}&\sum_{T_{i+1}^{(N)}\le t}\E\left[\Delta A_i^N\Delta L_{T_{i+1}^{(N)}}\Big|\mathcal{F}_{T_i^{(N)}}\right]\stackrel{p}{\longrightarrow} 0~,\\
&\sum_{T_{i+1}^{(N)}\le t}\E\left[\Delta A_i^N\Delta M_{T_{i+1}^{(N)}}\Big|\mathcal{F}_{T_i^{(N)}}\right]\stackrel{p}{\longrightarrow} 0~.\end{align*}
\end{lem}
\begin{proof}
Both relations are proven similarly and we restrict ourselves to the proof of the first one. The left-hand side equals
{\allowdisplaybreaks[2]{
\begin{align*}&
\sqrt{\frac{N}{T}}\sum_{T_{i+1}^{(N)}\le t}\E\left[\Delta L_{T_{i+1}}\left((L_{g_i}-L_{T_i})(M_{T_i}-M_{\lambda_i})+(M_{\gamma_i}-M_{T_i})(L_{T_i}-L_{l_i})
\right.\right.\\
&~~~~~~~~~~~~~\left.\left.+(L_{T_i}-L_{l_{i+1}})(M_{T_{i+1}}-M_{T_i})+(M_{T_i}-M_{\lambda_{i+1}})(L_{T_{i+1}}-L_{T_i})\right)\Big|\mathcal{F}_{{T_i}^{(N)}}\right]\\
&=\sqrt{\frac{N}{T}}\sum_{T_{i+1}^{(N)}\le t}\left(\E\left[\int_{T_i}^{g_i}(\sigma_s^X)^2ds\right](M_{T_i}-M_{\lambda_i})+\E\left[\int_{T_i}^{\gamma_i}\rho_s\sigma_s^X\sigma_s^Yds\right](L_{T_i}-L_{l_i})
\right.\\
&~~~~~~~~\left.+(L_{T_i}-L_{l_{i+1}})\E\left[\int_{T_i}^{T_{i+1}}\rho_s\sigma_s^X\sigma_s^Yds\right]+(M_{T_i}-M_{\lambda_{i+1}})\E\left[\int_{T_i}^{T_{i+1}}(\sigma_s^X)^2ds\right]\right)=\:\Gamma~.\end{align*}
$\Gamma$ is centred and using It\^{o} isometry the variance is shown to converge to zero:
\begin{align*}
&\frac{N}{T}\sum_{T_{i+1}^{(N)}\le t}\left(\left(\E\left[\int_{T_i}^{g_i}(\sigma_s^X)^2ds\right]\right)^2\E\left[\int_{\lambda_i}^{T_i}(\sigma_s^Y)^2ds\right]+\left(\E\left[\int_{T_i}^{\gamma_i}\rho_s\sigma_s^X\sigma_s^Yds\right]\right)^2\E\left[\int_{l_i}^{T_i}(\sigma_s^X)^2ds\right]
\right.\\
&~~~~~\left.+\E\left[\int_{l_{i+1}}^{T_i}(\sigma_s^X)^2ds\right]\left(\E\left[\int_{T_i}^{T_{i+1}}\rho_s\sigma_s^X\sigma_s^Yds\right]\right)^2+\E\left[\int_{\lambda_{i+1}}^{T_i}(\sigma_s^Y)^2ds\right]\left(\E\left[\int_{T_i}^{T_{i+1}}(\sigma_s^X)^2ds\right]\right)^2\right.\\
&~~~~~~~~~~\left. +2\E\left[\int_{l_{i+1}}^{T_i}\rho_s\sigma_s^X\sigma_s^Yds\right]\E\left[\int_{T_i}^{g_i}(\sigma_s^X)^2ds\right]\E\left[\int_{T_i}^{T_{i+1}}\rho_s\sigma_s^X\sigma_s^Yds\right]\right.\\
&~~~~~~~~~~\left.+2\E\left[\int_{\lambda_{i+1}}^{T_i}\rho_s\sigma_s^Y\sigma_s^Yds\right]\E\left[\int_{T_i}^{T_{i+1}}(\sigma_s^X)^2ds\right]\E\left[\int_{T_i}^{\gamma_i}\rho_s\sigma_s^X\sigma_s^Yds\right]\right)\\
&\le C N\delta_N^2=\KLEINO(1)~.\end{align*}}}
Once more we can conclude that the addends are uncorrelated since $\gamma_i=T_i\Rightarrow \lambda_{i+1}=T_i$ and $g_i=T_i\Rightarrow l_{i+1}=T_i$, respectively.
\end{proof}
Finally, we prove that the discrete covariation of our considered martingale with every bounded $\mathcal{F}_t$-martingale that is orthogonal to $L_t$ or $M_t$, converges to zero in probability. Hence, this lemma will complete the proof of Proposition \ref{HYasB}.
\begin{lem}\label{HYasC2}
Assume that $L_t^{\bot}$ and $M_t^{\bot}$ are bounded $\mathcal{F}_t$-martingales, with
$\left[ L,L^{\bot}\right] \equiv 0$ and $\left[ M,M^{\bot}\right]\equiv 0$, respectively. It holds true that
\begin{align*}&\sum_{T_{i+1}^{(N)}\le t}\E\left[\Delta A_i^N\Delta L^{\bot}_{T_{i+1}^{(N)}}\Big|\mathcal{F}_{T_i^{(N)}}\right]\stackrel{p}{\longrightarrow} 0~,\\
&\sum_{T_{i+1}^{(N)}\le t}\E\left[\Delta A_i^N\Delta M^{\bot}_{T_{i+1}^{(N)}}\Big|\mathcal{F}_{T_i^{(N)}}\right]\stackrel{p}{\longrightarrow} 0~.\end{align*}
\end{lem}
\begin{proof}
As in the preceding lemma, we only prove the first part of the result. The left-hand side of the first equation equals
{\allowdisplaybreaks[2]{
\begin{align*}&
\sqrt{\frac{N}{T}}\sum_{T_{i+1}^{(N)}\le t}\E\left[\Delta L^{\bot}_{T_{i+1}}\left((L_{g_i}-L_{T_i})(M_{T_i}-M_{\lambda_i})+(M_{\gamma_i}-M_{T_i})(L_{T_i}-L_{l_i})
\right.\right.\\
&~~~~~~~~~~~~~\left.\left.+(L_{T_i}-L_{l_{i+1}})(M_{T_{i+1}}-M_{T_i})+(M_{T_i}-M_{\lambda_{i+1}})(L_{T_{i+1}}-L_{T_i})\right)\Big|\mathcal{F}_{{T_i}^{(N)}}\right]\\
&=\sqrt{\frac{N}{T}}\sum_{T_{i+1}^{(N)}\le t}\left(\E\left[\int_{T_i}^{g_i}d\left[ L,L^{\bot}\right]_s\right](M_{T_i}-M_{\lambda_i})+\E\left[\int_{T_i}^{\gamma_i}d\left[ M,L^{\bot}\right]_s\right](L_{T_i}-L_{l_i})
\right.\\
&~~~~~~~~\left.+(L_{T_i}-L_{l_{i+1}})\E\left[\int_{T_i}^{T_{i+1}}d\left[ M,L^{\bot}\right]_s\right]+(M_{T_i}-M_{\lambda_{i+1}})\E\left[\int_{T_i}^{T_{i+1}}d\left[ L,L^{\bot}\right]_s\right]\right)\\
&=\sqrt{\frac{N}{T}}\sum_{T_{i+1}^{(N)}\le t}\E\left[\int_{T_i}^{T_{i+1}}d\left[ M,L^{\bot}\right]_s\right](L_{T_i}-L_{l_{i+1}})+\E\left[\int_{T_i}^{\gamma_i}d\left[ M,L^{\bot}\right]_s\right](L_{T_i}-L_{l_i})~.\end{align*}}}
This term is centred and the has the variance
{\allowdisplaybreaks[2]{
\begin{align*}&\frac{N}{T}\sum_{T_{i+1}^{(N)}\le t}\E\left[\left(\int_{T_i}^{T_{i+1}}d\left[ M,L^{\bot}\right]_s\right)^2\right]
\E\left[\left(L_{T_i}-L_{l_{i+1}}\right)^2\right]
+\E\left[\left(\int_{T_i}^{\gamma_i}d\left[ M,L^{\bot}\right]_s\right)^2\right]\\&~~\times\E\left[\left(L_{T_i}-L_{l_{i}}\right)^2\right]+\E\left[(L_{T_i}-L_{L_{i+1}})^2\right]\E\left[\int_{T_i}^{T_{i+1}}d\left[ M,L^{\bot}\right]_s\right]\E\left[\int_{T_i}^{\gamma_i}d\left[ M,L^{\bot}\right]_s\right]\\
&{\underset{\text{It\^{o} isometry}}{=}}\frac{N}{T}\sum_{T_{i+1}^{(N)}\le t}\E\left[\left(\int_{T_i}^{T_{i+1}}d\left[ M,L^{\bot}\right]_s\right)^2\right]\E\left[\int_{l_{i+1}}^{T_i}(\sigma_s^X)^2ds\right]\\
&~~~~~~~~~~~~~~~~~~~~~~~~~+\E\left[\left(\int_{T_i}^{\gamma_i}d\left[ M,L^{\bot}\right]_s\right)^2\right]\E\left[\int_{l_i}^{T_i}(\sigma_s^X)^2ds\right]\\
&~~~~~~~~~~~~~~~~~~~~~~~~~+\E\left[\int_{l_{i+1}}^{T_i}(\sigma_s^X)^2ds\right]\E\left[\int_{T_i}^{T_{i+1}}d\left[ M,L^{\bot}\right]_s\right]\E\left[\int_{T_i}^{\gamma_i}d\left[ M,L^{\bot}\right]_s\right]=\KLEINO(1)~.\end{align*}}}
Thus, the covariations converge to zero in probability.
\end{proof}
The Lemma completes the proof of Proposition \ref{HYasB}.
\end{proof}
The mixed normal limit in Proposition \ref{HYas} can be obtained as the marginal distribution of $\mathcal{A}_T^N$ in $t=T$. 
\end{proof}
Proposition \ref{HYas} for the error of the approximation by the discretization error of the closest synchronous approximation \eqref{As} and the stable limit theorem for this synchronous discretization error \eqref{D} given in Proposition \ref{HYdis} suffice to imply Theorem \ref{HYclt}.
That is because the multivariate stable convergence Theorem \ref{jacodmult} applies to the vector of the two uncorrelated terms and since the covariations converge to zero , the stable convergence to the mixed Gaussian limit with the sum of the two asymptotic variances is concluded.$\hfill\Box$\\ \noindent

\section{Proof of Proposition \ref{avarest}\label{sec:8}}
The proof will be divided into three parts in that the sum of squared products, products of consecutive increments and the histogram estimator are considered, respectively. Denote $X_j^+=X_{g_j}-X_{T_j}$, $X_j^-=X_{T_{j-1}}-X_{l_j}$ and $X_j^S=X_{T_j}-X_{T_{j-1}}~,~j=1,\ldots,N$. In the first step it is proved that
$$N\hspace*{-.05cm}\sum_{j=1}^{N-1}\left(X_{g_j}-X_{l_j}\right)^2\left(Y_{\gamma_j}-Y_{\lambda_j}\right)^2\hspace*{-.1cm}\stackrel{p}{\longrightarrow}\hspace*{-.1cm}T\hspace*{-.1cm}\int_0^T\hspace*{-.15cm} G^{\prime}(t)\hspace*{-.05cm}\left(\sigma_t^X\sigma_t^Y\right)^2\hspace*{-.05cm}\left(2\rho_t^2+1\right)dt+T\hspace*{-.1cm}\int_0^T\hspace*{-.15cm}F^{\prime}(t)\hspace*{-.05cm}\left(\sigma_t^X\sigma_t^Y\right)^2\hspace*{-.05cm}dt~.$$
\begin{align*}N\hspace*{-.1cm}\sum_{j=1}^{N-1}\hspace*{-.05cm}\left(X_j^++X_j^S+X_j^-\right)^2\left(Y_j^++Y_j^S+Y_j^-\right)^2&\hspace*{-.1cm}=\hspace*{-.05cm}N\hspace*{-.05cm}\sum_{j=1}^{N-1}\hspace*{-.05cm}\left((X_j^+)^2(Y_j^S+Y_j^-)^2\hspace*{-.05cm}+\hspace*{-.05cm}(Y_j^+)^2(X_j^S+X_j^-)^2\right.\\
&\left. \hspace*{-.15cm}+(X_j^-)^2(Y_j^S)^2\hspace*{-.05cm}+\hspace*{-.05cm}(Y_j^-)^2(X_j^S)^2\hspace*{-.05cm}+\hspace*{-.05cm}(X_j^SY_j^S)^2\right)\hspace*{-.05cm}+\hspace*{-.1cm}\KLEINO_p(1)\end{align*}
All centred addends have a variance tending to zero as $N\rightarrow\infty$ and converge to zero in probability. The sum of the first four addends times the factor $N/T$ has been proved to converge in probability to $\int_0^TF^{\prime}(t)(\sigma_t^X\sigma_t^Y)^2dt$ in Lemma \ref{HYasV} where this term has appeared in the sequence of conditional variances of the error due to non-synchronicity.\\
Hence, it remains to prove that $N\sum_{j=1}^{N-1}(X_j^SY_j^S)^2\stackrel{p}{\rightarrow}T\int_0^T(2\rho_t^2+1)(\sigma_t^X\sigma_t^Y)^2G^{\prime}(t)dt$.
For this purpose recall the notation from the proof of Proposition \ref{HYdis}. With $L_t=\int_0^t\sigma_s^XdW_s^X~,~M_t=\int_0^t\sigma_s^YdW_s^Y$, $L_i=L_{T_i},M_i=M_{T_i}$, we can write the term
\begin{align*}N\sum_{j=1}^{N-1}\left((L-L_{i-1})_{T_i}(M-M_{i-1})_{T_i}\right)^2&=N\sum_{i=1}^{N-1}\left(2\int_0^{T_i}(L-L_{i-1})_t(M-M_{i-1})_t^2d(L-L_{i-1})_t\right.\\
&\left. ~+2\int_0^{T_i}(L-L_{i-1})^2_t(M-M_{i-1})_td(M-M_{i-1})_t\right.\\
&\hspace*{-.4cm}\left. ~+4\int_0^{T_i}(L-L_{i-1})_t(M-M_{i-1})_td\left[ M-M_{i-1},L-L_{i-1}\right]_t\right.\\
&\hspace*{-2cm}\left.+\int_0^{T_i}(M-M_{i-1})_t^2 d\left[ M-M_{i-1}\right]_t+\int_0^{T_i}(L-L_{i-1})_t^2d\left[ L-L_{i-1}\right]_t\right)~,\end{align*}
where we have applied It\^{o}'s formula.
The sum of the first two addends converges to zero in probability since it is centred and the variance converges to zero. Since
$$\int_0^{T_i}(L-L_{i-1})_t(M-M_{i-1})_td\left[ M-M_{i-1},L-L_{i-1}\right]_t=\int_{T_{i-1}}^{T_i}(L-L_{i-1})_t(M-M_{i-1})_td\left[ M,L\right]_t~,$$
the sum of the third addends has been considered in the proof of Proposition \ref{HYdisB} as part of the quadratic variation of the discretization error of the closest synchronous approximation and converges in probability to $2\,T\int_0^TG^{\prime}(t)(\rho_t\sigma_t^X\sigma_t^Y)^2dt$. The remaining sum of the fourth addends is also similar to the other part of the quadratic variation in the proof of Proposition \ref{HYdisB}. An analogous approximation and integration by parts yields
\begin{align*}&\int_{T_{i-1}}^{T_i}(M-M_{i-1})_t^2d\left[ M\right]_t+\int_{T_{i-1}}^{T_i}(L-L_{i-1})_t^2d\left[ L\right]_t\\ &~~=\int_{T_{i-1}}^{T_i}\left[ M-M_{i-1}\right]_td\left[ M\right]_t+\int_{T_{i-1}}^{T_i}\left[ L-L_{i-1}\right]_td\left[ L\right]_t+\KLEINO_p(1)\\
\\&~~=\int_{T_{i-1}}^{T_i}d(\left[ L-L_{i-1}\right]_t\left[ M-M_{i-1}\right]_t)+\KLEINO_p(1)\end{align*}
and the convergence of the above given term to $T\int_0^T(\sigma_t^X\sigma_t^Y)^2G^{\prime}(t)dt$.\\
In the second part of the proof we are concerned with the term
\begin{align*}2N\sum_{j=1}^{N-1}(X_j^++X_j^S+X_j^-)(X_{j+1}^++X_{j+1}^S+X_{j+1}^-)(Y_j^++Y_j^S+Y_j^-)(Y_{j+1}^++Y_{j+1}^S+Y_{j+1}^-)\\
=2N\sum_{j=1}^{N-1}\left(X_j^SY_j^SX_{j+1}^SY_{j+1}^S+X_j^+Y_j^SX_{j+1}^-Y_{j+1}^S+Y_j^+X_j^SY_{j+1}^-X_{j+1}^S\right)+\KLEINO_p(1)~.\end{align*}
The sum incorporating all centred addends converges to zero in probability. The last two addends capture the only dependence between consecutive addends in the error due to non-synchronicity \eqref{As}, namely when next-tick interpolations and previous-tick interpolations at the same $T_i,i=1,\ldots,N$ are included. Those have appeared in the proof of Lemma \ref{HYasV} and have been proved to converge to $T\int_0^T2H^{\prime}(t)(\sigma_t^X\sigma_t^Y)^2dt$ in probability. That $2N\sum(X_j^SY_j^SX_{j+1}^SY_{j+1}^S)\stackrel{p}{\rightarrow}2\int_0^TG^{\prime}(t)(\rho_t\sigma_t^X\sigma_t^Y)^2dt$ follows with 
the methodology from \cite{bibinger2} and Lemma 1 from \cite{zhang} using the concept of a time-change in the asymptotic quadratic variation of refresh times such that $\sum_i(\Delta T_i-T/N)^2=\KLEINO(N^{-1})$ holds true. Using the mean value theorem and $(\Delta T_j)^2-\Delta T_j\Delta T_{j+1}=\Delta T_j(\Delta T_j-T/N)+\Delta T_j(T/N-\Delta T_{j+1})$ together with the Cauchy-Schwarz inequality
$$N\left|\sum_{j=1}^{N-1}\left(\Delta T_j\left(\Delta T_j-\frac{T}{N}\right)\right)\right|\le N\sqrt{\sum_{j=1}^{N-1}(\Delta T_j)^2}\sqrt{\sum_{j=1}^{N-1}\left(\Delta T_j-\frac{T}{N}\right)^2}$$
yields the result.\\
The Hayashi-Yoshida estimators on the bins in the histogram-based estimator \eqref{hist} fulfill
$$\widehat{\Delta\left[ X,Y\right]}_{G_j^N}^{(HY)}=\int_{G_{j-1}^N}^{G_j^N}\rho_t\sigma_t^X\sigma_t^Ydt+\mathcal{O}_p\left(K_N^{\nicefrac{1}{2}}N^{-\nicefrac{1}{2}}\right)$$
so that the estimation error of the sum is of order $K_N^{\nicefrac{3}{4}}N^{-\nicefrac{1}{2}}$ in probability and for $K_N\rightarrow\infty\,,\,N\rightarrow\infty\,,\,K_NN^{-\nicefrac{2}{3}}\rightarrow 0$\sectionmark{Asymptotic variance estimation} consistency holds and we conclude consistency of the estimator of the asymptotic variance.$\hfill\Box$\\ \noindent

\bibliographystyle{model1b-num-names}
\bibliography{literatur}

\end{document}